\pdfoutput=1
\RequirePackage{ifpdf}
\ifpdf 
\documentclass[pdftex]{sigma}
\else
\documentclass{sigma}
\fi

\numberwithin{equation}{section}

\newtheorem{Theorem}{Theorem}[section]
\newtheorem{Corollary}[Theorem]{Corollary}

\newtheorem{Proposition}[Theorem]{Proposition}
{ \theoremstyle{definition}

\newtheorem{Remark}[Theorem]{Remark} }

\def\al{\alpha}
\def\be{\beta}
\def\de{\delta}
\def\la{\lambda}
\def\si{\sigma}
\def\ep{\varepsilon}
\def\vp{\varphi}
\def\om{\omega}
\def\Ga{\Gamma}
\def\De{\Delta}
\def\R{\mathbb{R}}
\def\N{\mathbb{N}}
\def\cD{\mathcal D}
\def\cV{\mathcal V}
\def\Dom{\mathrm{Dom}}
\newcommand{\qbinom}[3]{{\genfrac{[}{]}{0pt}{}{#1}{#2}}_{#3}}
\newcommand{\rFs}[5]{\,_{#1}F_{#2} \left( \genfrac{.}{.}{0pt}{}{#3}{#4};#5 \right)}
\newcommand{\rphis}[5]{\,_{#1}\varphi_{#2} \left( \genfrac{.}{.}{0pt}{}{#3}{#4} ;#5 \right)}

\begin{document}

\allowdisplaybreaks

\newcommand{\arXivNumber}{1802.09190}

\renewcommand{\thefootnote}{}

\renewcommand{\PaperNumber}{072}

\FirstPageHeading

\ShortArticleName{Generalized Burchnall-Type Identities for Orthogonal Polynomials and Expansions}

\ArticleName{Generalized Burchnall-Type Identities\\ for Orthogonal Polynomials and Expansions\footnote{This paper is a~contribution to the Special Issue on Orthogonal Polynomials, Special Functions and Applications (OPSFA14). The full collection is available at \href{https://www.emis.de/journals/SIGMA/OPSFA2017.html}{https://www.emis.de/journals/SIGMA/OPSFA2017.html}}}

\Author{Mourad E.H.~ISMAIL~$^\dag$, Erik KOELINK~$^\ddag$ and Pablo ROM\'AN~$^\S$}

\AuthorNameForHeading{M.E.H.~Ismail, E.~Koelink and P.~Rom\'an}

\Address{$^\dag$~University of Central Florida, Orlando, Florida 32816, USA}
\EmailD{\href{mailto:ismail@math.ucf.edu}{ismail@math.ucf.edu}}

\Address{$^\ddag$~IMAPP, Radboud Universiteit, PO Box 9010, 6500GL Nijmegen, The Netherlands}
\EmailD{\href{mailto:e.koelink@math.ru.nl}{e.koelink@math.ru.nl}}

\Address{$^\S$~CIEM, FaMAF, Universidad Nacional de C\'ordoba,\\
\hphantom{$^\S$}~Medina Allende s/n Ciudad Universitaria, C\'ordoba, Argentina}
\EmailD{\href{mailto:roman@famaf.unc.edu.ar}{roman@famaf.unc.edu.ar}}

\ArticleDates{Received February 27, 2018, in final form July 11, 2018; Published online July 17, 2018}

\Abstract{Burchnall's method to invert the Feldheim--Watson linearization formula for the Hermite polynomials is extended to all polynomial families in the Askey-scheme and its $q$-analogue. The resulting expansion formulas are made explicit for several families corresponding to measures with infinite support, including the Wilson and Askey--Wilson polynomials. An integrated version gives the possibility to give alternate expression for orthogonal polynomials with respect to a modified weight. This gives expansions for polynomials, such as Hermite, Laguerre, Meixner, Charlier, Meixner--Pollaczek and big $q$-Jacobi polynomials and big $q$-Laguerre polynomials. We show that one can find expansions for the orthogonal polynomials corresponding to the Toda-modification of the weight for the classical polynomials that correspond to known explicit solutions for the Toda lattice, i.e., for Hermite, Laguerre, Charlier, Meixner, Meixner--Pollaczek and Krawtchouk polynomials.}


\Classification{33C45; 33D45; 42C05; 37K10}

\renewcommand{\thefootnote}{\arabic{footnote}}
\setcounter{footnote}{0}

\section{Introduction}\label{sec:intro}

In 1941 J.L.~Burchnall (1892--1975) wrote a short paper \cite{Burc} in which he developed a method to find the inverse to the 1938 Feldheim--Watson formula
\begin{gather}\label{eq:intro-FeldheimWatson}
H_m(x)H_n(x) = \sum_{r=0}^{m\wedge n} \binom{n}{r} \binom{m}{r}2^r r! H_{m+n-2r}(x),
\end{gather}
linearizing a product of two Hermite polynomials. Here we use standard notation for the Hermite polynomials $H_m$, and see \cite[Section~6.1]{AndrAR}, \cite[Lecture~5]{Aske-SIAM}, \cite[Section~4.6]{Isma}, The composition of Burchnall's formula \cite[equation~(5)]{Burc}
\begin{gather}\label{eq:exCOP-Hermite-Burchnall2}
H_{n+m}(x) = \sum_{r=0}^{n\wedge m} \binom{n}{r} \binom{m}{r} (-2)^{r} r! H_{n-r}(x) H_{m-r}(x),
\end{gather}
which we rederive in Section~\ref{ssec:Hermite}, and the Feldheim--Watson formula \eqref{eq:intro-FeldheimWatson} is equivalent to the finite Chu--Vandermonde sum. In fact, Nielsen \cite[p.~33, equations~(5) and (6)]{Niel} derived the Feldheim--Watson and the Burchnall formulas already in 1918. It seems that Burchnall was not aware of Nielsen's paper. Nielsen used the recurrence relations, whereas Burchnall used an operational approach. In fact in his series of memoirs in the 1890s, L.J.~Rogers calculated the linearization coefficients for the continuous $q$-ultraspherical polynomials and the limiting case of the continuous $q$-Hermite polynomials, see, e.g., \cite[Section~8.5]{GaspR} and references given there. So~\eqref{eq:intro-FeldheimWatson} is a $q \to 1$ limiting case of one of Rogers's results. Rogers used his results to prove the Rogers--Ramanujan identities, see, e.g., \cite[Section~8.10]{GaspR}.

Actually, Burchnall gives an operational formula for $\big(\frac{d}{dx}-2x\big)^n$ in terms of the $k$-th derivatives. Explicitly, see \cite[equation~(3)]{Burc},
\begin{gather}\label{eq:exCOP-Hermite-Burchnall}
\left( \left(\frac{d}{dx}- 2x\right)^n f\right)(x)=\sum_{k=0}^n(-1)^{n-k} \binom{n}{k} H_{n-k}(x) f^{(k)}(x),
\end{gather}
for $f$ sufficiently differentiable. The important ingredient in Burchnall's derivation are the raising and lowering operators for the Hermite polynomials. In Section~\ref{sec:general-setup} we show that Burchnall's method extends to all polynomials in the Askey-scheme and its $q$-analogue, see, e.g., \cite{KoekLS,KoekS}, using their raising and lowering operators. We generalize~\eqref{eq:exCOP-Hermite-Burchnall} in Theorem~\ref{lem:gs-expansion}. We give many worked out cases of the analogue of Burchnall's identity up to the Wilson and the Askey--Wilson polynomials. We restrict ourselves to families of orthogonal polynomials with infinite support, although it is clear that the method works for finite discrete orthogonal polynomials as well, up to minor modifications.

It should be noted that Carlitz \cite[equation~(4)]{Carl} extends Burchnall's operational formula \eqref{eq:exCOP-Hermite-Burchnall} to the case of Laguerre polynomials which we rederive in \eqref{eq:Laguerre-lemmagsexpansion}. The resulting analogue of \eqref{eq:exCOP-Hermite-Burchnall2} is
\begin{gather}\label{eq:Laguerre-corlemmagsexpansion}
\frac{(m+1)_n }{n!} L_{m+n}^{(\nu)}(x) = \sum_{k=0}^n\frac{(-x)^k}{k!}L_{n-k}^{(\nu+k)}(x) L_{m-k}^{(\nu+n+k)}(x),
\end{gather}
which Carlitz \cite{Carl} proves by induction on $n$. Gould and Hopper \cite{GoulH} give a joint extension of these operational formulas to the Gould--Hopper polynomials. Al-Salam \cite{AlSa} shows that the operational formulas of Carlitz and Gould--Hopper are equivalent. Moreover, Singh \cite{Sing} has determined the extension to the Jacobi case, which we give in Section~\ref{ssec:Jacobi} for completeness. The extensions of Burchnall's results to the Laguerre and Jacobi case fit in the classical orthogonal polynomials as a part of the Askey scheme. This part is characterized by having the derivative as the lowering operator. The extension to Zernike polynomials is given in~\cite{AharAFG}, where more references to the literature are given.

In this paper we make the operational formulation \eqref{eq:exCOP-Hermite-Burchnall} and the analogue of the expansion formulas \eqref{eq:exCOP-Hermite-Burchnall2}, \eqref{eq:Laguerre-corlemmagsexpansion} explicit for several classes in the Askey scheme and its $q$-analogue. For the classical orthogonal polynomials~-- Hermite, Laguerre, Jacobi~-- this is in Section~\ref{sec:examples-classicalOPS}, and we show that the integrated formulas in case of the Hermite and Laguerre polynomials are essentially given as a change of coordinates. As is clear from the above, several of these results for classical orthogonal polynomials can be traced back in the literature, see \cite{AharAFG, AlSa,Burc,Carl,GoulH,Niel,Sing} and references given there. We show how to generalize these operational formulas and expansion identities to all of the families of orthogonal polynomials in the Askey scheme and its $q$-analogue. Replacing the derivative by the backward shift operator as the lowering operator, we get two versions of the operational identities and the expansion formulas for the corresponding classes of orthogonal polynomials: Meixner and Charlier polynomials, see Section~\ref{sec:backwardshift}. In these cases there are two inequivalent versions of the expansion formulas of type \eqref{eq:exCOP-Hermite-Burchnall2}, \eqref{eq:Laguerre-corlemmagsexpansion}. Using the divided difference operator $\frac{\delta}{\delta x}$ we give the explicit formulas for the Meixner--Pollaczek polynomials in Section~\ref{sec:div-difference}. For the Wilson polynomials related to the divided difference operator $\frac{\delta}{\delta x^2}$ the results are stated in Section~\ref{sec:examples-Wilson}. For the $q$-analogue we study the big $q$-Jacobi polynomials and the Askey--Wilson polynomials. For the big $q$-Jacobi polynomials we again get two versions of the operational identity and the expansion formula given in Section~\ref{sec:examples-qAskey-qderivative}. By switching to the Askey--Wilson $q$-difference operator as a lowering operator, we may include the Askey--Wilson polynomials. This is in Section~\ref{sec:AskeyWilsoncase}, with emphasis on the special case of the continuous $q$-Hermite polynomials.

Using the generalized Burchnall identity and the raising and lowering operators being each others adjoint on suitable weighted $L^2$-spaces, we obtain an integrated version in Corollary~\ref{cor:lem:gs-expansion2}. In various cases we can use the integrated version to find orthogonal polynomials with respect to a slightly modified weight, where the modification is governed by an eigenfunction to certain operators occurring in the generalized Leibniz rule. We use this to find the orthogonal polynomials with respect to the orthogonality measure modified by an exponential or a $q$-exponential function. Often it is possible to recognize the modified measure explicitly in typically the same class. In this way we obtain expressions for orthogonal polynomials in terms of orthogonal polynomials with different parameters. In particular, for the modification by $e^{-xt}$ we can understand the orthogonal polynomials for the measure $e^{-xt}d\mu(x)$, where the measure $\mu$ corresponds to the orthogonality measure for a family of orthogonal polynomials in the Askey-scheme. By the work of Flaschka, Moser and others this is related to the Lax pair for the Toda lattice, see, e.g., \cite[Section~4.6]{BabeBT}, \cite[Section~2]{BasoCE}, \cite[Section~2.8]{Isma}, \cite[Section~2]{Zhed-R}. The recurrence coefficients for the monic orthogonal polynomials $x p_n(x;t) = p_{n+1}(x;t) + b_n(t) p_n(x;t) + c_n(t)p_{n-1}(x;t)$ for the modified weight satisfy the Toda lattice equations
\begin{gather}\label{eq:intro-Toda}
\dot{c}_n(t) = c_n(t)\bigl( b_{n-1}(t) - b_n(t)\bigr),\qquad \dot{b}_n(t) = c_n(t) - c_{n+1}(t).
\end{gather}
This approach works for the Hermite, Laguerre, Krawtchouk, Meixner, Charlier, and Meixner--Pollaczek polynomials, and this is listed in Proposition~\ref{prop:explicit-Toda-solutions}. Proposition~\ref{prop:explicit-Toda-solutions} is well-known, and Zhedanov \cite{Zhed} derives the result from the requirement that the Lax pair $(L,M)$ and the time-derivative $\dot{L}$ close up to a $3$-dimensional Lie algebra. In this context, Proposition~\ref{prop:explicit-Toda-solutions} follows from Lie algebra representations, see, e.g., \cite[Section~4.6]{BabeBT} for the link of the Toda lattice to simple Lie algebras, and its relation to orthogonal polynomials, see, e.g., \cite{KoelVdJ} for the link to the corresponding polynomials.

For the $q$-exponential functions $e_q$ and $E_q$ we find expansions of the big $q$-Jacobi polynomials in terms of big $q$-Laguerre polynomials. We also find the inverse formulas in this way. However, there seems no integrable system associated to these expansions.

This paper deals with the scalar case. In a companion paper~\cite{IsmaKR} we also use Burchnall type identities for matrix-valued orthogonal polynomials in order to give a non-trival solution to the non-abelian Toda lattice analogue of~\eqref{eq:intro-Toda} as introduced by Gekhtman~\cite{Gekh}.

The contents of the paper are related to the role of the analogue of differentiation in the Askey-scheme and its $q$-analogue. First, in Section~\ref{sec:general-setup} we write down the general set-up as motivated by Burchnall's paper \cite{Burc} and extending it. In Section~\ref{sec:examples-classicalOPS} we look at those polynomials that correspond to the (ordinary) derivative, which are the classical polynomials of Hermite, Laguerre and Jacobi. We recover many of the known results in this way. By specializing we also recover expansions for Laguerre and Hermite polynomials, which can be considered as special cases of more general convolution identities \cite{KoelVdJ}. These expansions are related to explicit solutions of the Toda lattice \eqref{eq:intro-Toda}. In Section~\ref{sec:backwardshift} we consider the shift operator as the analogue of the derivative, and we consider the families of Meixner and Charlier polynomials. In Section~\ref{sec:div-difference} we consider the Meixner--Pollaczek polynomials, and in the following Section~\ref{sec:examples-Wilson} we consider the Wilson polynomials. In Section~\ref{sec:Todalattice} we summarize the explicit solutions of the Toda lattice arising in this way. In the remaining sections we focus on the $q$-Askey scheme. In Section~\ref{sec:examples-qAskey-qderivative} we consider the $q$-Hahn scheme with the big $q$-Jacobi polynomials on top. Here we give several expansions for polynomial families in the $q$-Hahn scheme. In Section~\ref{sec:AskeyWilsoncase} we discuss the top-level of the Askey--Wilson polynomials.

We use standard notation for Pochhammer symbols, $q$-shifted factorials, hypergeometric and basic hypergeometric series, etc.\ as in, e.g., \cite{AndrAR, GaspR,Isma, KoekLS, KoekS, Rain}. We take $0<q<1$ for the base $q$, and we follow the convention that, e.g., $(x/ab;q)_\infty=\big(xa^{-1}b^{-1};q\big)_\infty$.

\section{Generalized Burchnall-type identities}\label{sec:general-setup}

In this section we show how the Burchnall identity for the Hermite polynomials as well as Burchnall's approach \cite{Burc} can be generalized to families of orthogonal polynomials for which raising and lowering operators exist. In particular we give an expansion for a class of functions in terms of the polynomials involved. Note that all the polynomial families in the Askey-scheme and its $q$-analogue with infinite support fulfill the assumptions in this section up to Corollary~\ref{cor:lem:gs-expansion2}, so that in particular Theorem~\ref{lem:gs-expansion} and Corollaries~\ref{cor:lem:gs-expansion},~\ref{cor:lem:gs-expansion2} are valid in all cases studied in this paper.

Let the general orthogonal polynomials satisfy
\begin{gather*} 
\int_\R p^{(\nu)}_n(x) p^{(\nu)}_m(x) d\mu^{(\nu)}(x) = \de_{m,n} h^{(\nu)}_n.
\end{gather*}
Here $\mu^{(\nu)}$ is a positive Borel measure on the real line for which all moments exist. The parame\-ter~$\nu$ is contained in some (multivariable) parameter set $\cV\subset \R^n$. We write a combination of parameters additively for the Askey scheme, but it is more convenient to write it multiplicatively for the $q$-Askey scheme. We assume that the measure $\mu^{(\nu)}$ has infinite support, but similar results for finite discrete orthogonal polynomials can be obtained up to some modifications.

Next we assume that for two elements $\nu, \nu+\si\in \cV$ of the parameter space we have (densely defined) raising operators $R_\nu \colon L^2\big(\mu^{(\nu+\si)}\big) \to L^2\big(\mu^{(\nu)}\big)$. Here $\si$ is a fixed element in $\R^n$, so that $\nu\in \cV$ implies that $\nu+\N \si \subset \cV$. In particular, we assume $R_\nu$ maps polynomials of degree $n$ to polynomials of degree $n+1$. The (densely defined) adjoint $L_\nu = (R_\nu)^\ast \colon L^2\big(\mu^{(\nu)}\big) \to L^2\big(\mu^{(\nu+\si)}\big)$ is assumed to be a~lowering operator, i.e., mapping polynomials of degree $n$ to polynomials of degree $n-1$. So in particular,
\begin{gather}\label{eq:gs-LRadjoint}
\langle R_\nu f, g \rangle_{L^2(\mu^{(\nu)})} =\langle f, L_\nu g \rangle_{L^2(\mu^{(\nu+\si)})}
\end{gather}
for all $f\in \Dom(R_\nu)\subset L^2\big(\mu^{(\nu+\si)}\big)$ and all $g\in \Dom(L_\nu)\subset L^2\big(\mu^{(\nu)}\big)$, and in particular we assume that the polynomials are contained in the domains of~$R_\nu$ and~$L_\nu$.

It follows that for $\nu\in \cV$ we can take
\begin{gather}\label{eq:gs-defOPs}
p_n^{(\nu)}(x) = \big(R_{\nu}R_{\nu+\si} \cdots R_{\nu+(n-2)\si} R_{\nu+(n-1)\si} \mathbf{1}\big)(x),
\end{gather}
where $\mathbf{1}(x)=1$ is the constant function. Note that in general $[R_{\nu+k\si}, R_{\nu+l\si}]\not=0$, so the order in~\eqref{eq:gs-defOPs} is relevant. Indeed, iterating~\eqref{eq:gs-LRadjoint} we find for polynomials~$f$,~$g$ that
\begin{gather}
\langle R_{\nu}R_{\nu+\si} \cdots R_{\nu+(n-2)\si} R_{\nu+(n-1)\si} f, g \rangle_{L^2(\mu^{(\nu)})} \nonumber\\
\qquad{}{} =\langle f, L_{\nu+(n-1)\si} L_{\nu+(n-2)\si} \cdots L_{\nu+\si} L_{\nu} g \rangle_{L^2(\mu^{(\nu+n\si)})}.\label{eq:gs-LRadjointiterated}
\end{gather}
Take $f=\mathbf{1}$, $g(x)=x^k$, $k<n$, so that $L_{\nu+(n-1)\si} L_{\nu+(n-2)\si} \cdots L_{\nu+\si} L_{\nu} g=0$ as a polynomial of negative degree $k-n$.

Taking $g(x)=p^{(\nu)}_n(x)$, $f(x)=x^k$, $k<n-1$, in \eqref{eq:gs-LRadjoint} we see that \smash{$\big\langle L_\nu p^{(\nu)}_n, x^k\big\rangle_{L^2(\mu^{(\nu+\si)})} =0$}, hence $L_\nu p^{(\nu)}_n = \ell^{(\nu)}_n p^{(\nu+\si)}_{n-1}$ for some constant $\ell^{(\nu)}_n$ which follows by~\eqref{eq:gs-defOPs} as $L_\nu R_\nu x^{n-1} = \ell^{(\nu)}_n x^{n-1} + \text{l.o.t.}$, where $x^n$, respectively $x^{n-1}$, can be replaced by some other suitable polynomials of degree~$n$, respectively~$n-1$, with leading coefficient~$1$. Here, and elsewhere $\text{l.o.t.}$ means `lower order terms'.

Next we assume that the raising operators have a specific form. So we assume that there exists a function $w_\nu$ for $\nu \in \cV$ so that $w_\nu>0$ for $\mu^{(\nu)}$-a.e.\ and that
\begin{gather*} 
R_\nu = M_{\nu}^{-1} \circ \partial \circ M_{\nu+\si}, \qquad (M_\nu f)(x) = (M_{w_\nu} f)(x) = w_\nu(x)f(x),
\end{gather*}
and $\partial$ is an operator independent of $\nu$. In particular, we assume that $w_\nu$ is in the domain of $\partial$. Typically, $\mu^{(\nu)}$ is absolutely continuous with respect to Lebesgue measure or a counting measure (independent of $\nu$), and the Radon--Nikodym derivative is $w^{(\nu)}$. So
\begin{gather*} 
R_{\nu}R_{\nu+\si} \cdots R_{\nu+(n-2)\si} R_{\nu+(n-1)\si} =M_{\nu}^{-1} \circ \partial^n \circ M_{\nu+n\si}.
\end{gather*}
Typically, $\partial$ is a lowering operator.

Moreover, we assume of the existence of a Leibniz rule for the operator $\partial$ of the form
\begin{gather}\label{eq:gs-Leibniz}
\bigl(\partial^n (fg)\bigr) (x) = \sum_{k=0}^n \al^n_k \big(\eta^k\partial^{n-k}f\big)(x) (T_{k,n}g)(x),
\end{gather}
where $\al^n_k$ are constants, like binomial coefficients, $\eta$ is a fixed operator (e.g., in order to accommodate the Askey--Wilson $q$-difference operator), and~$\eta$ can also be the identity. Here~$f$ and~$g$ are such that all expressions are well-defined. In general, we assume that~$\eta$ is an invertible homomorphism, so $\eta(fg) = \eta(f)\eta(g)$. The $T_{k,n}$ are analogues of differential operators, and for the general case we have to allow for $n$-dependence.

For a suitable function $f$, such as a holomorphic function on a sufficiently large domain, an entire function or a polynomial, we have
\begin{gather*}
\big(R_{\nu}R_{\nu+\si} \cdots R_{\nu+(n-2)\si} R_{\nu+(n-1)\si} f\big)(x) = \big(M_{\nu}^{-1} \circ \partial^n \circ M_{\nu+n\si} f\big)(x) \\
\qquad{} = \frac{1}{w_{\nu}(x)} \partial^n\bigl(w_{\nu+n\si} f\bigr) (x)= \frac{1}{w_{\nu}(x)} \sum_{k=0}^n \al^n_k \big(\eta^k\partial^{n-k}w_{\nu+n\si} \big)(x) (T_{k,n}f)(x),
\end{gather*}
and rewriting
\begin{gather*}
\big(\eta^k\partial^{n-k}w_{\nu+n\si} \big)(x) =\eta^k (w_{\nu+k\si})(x) \eta^k \big(p_{n-k}^{(\nu+k\si)}\big)(x)
\end{gather*}
proves Theorem~\ref{lem:gs-expansion}.

\begin{Theorem}\label{lem:gs-expansion}
For $f$ a suitable function, such as, e.g., a polynomial or a holomorphic function, we have
\begin{gather*}
\big(R_{\nu}R_{\nu+\si} \cdots R_{\nu+(n-2)\si} R_{\nu+(n-1)\si} f\big)(x) =\sum_{k=0}^n \al^n_k \frac{\eta^k (w_{\nu+k\si})(x)}{w_{\nu}(x)}\eta^k \big(p_{n-k}^{(\nu+k\si)}\big)(x)(T_{k,n}f)(x)
\end{gather*}
using the notation and assumptions as above.
\end{Theorem}

Theorem \ref{lem:gs-expansion} is the extension of Burchnall's formula \cite[equation~(3)]{Burc} for Hermite polynomials, see~\eqref{eq:exCOP-Hermite-Burchnall}.

Take $f(x) = (R_{\nu+n\si} \cdots R_{\nu+(n+m-1)\si} \mathbf{1})(x)= p^{(\nu+n\si)}_{m}(x)$ using~\eqref{eq:gs-defOPs}, then we find from Theo\-rem~\ref{lem:gs-expansion} and~\eqref{eq:gs-defOPs} the following polynomial identity, again using the notation of this section.

\begin{Corollary}\label{cor:lem:gs-expansion} With the notation and assumptions as above
\begin{gather*}
p^{(\nu)}_{n+m}(x) =\sum_{k=0}^n \al^n_k \frac{\eta^k (w_{\nu+k\si})(x)}{w_{\nu}(x)}\eta^k \big(p_{n-k}^{(\nu+k\si)}\big)(x)\big(T_{k,n}p^{(\nu+n\si)}_{m}\big)(x).
\end{gather*}
\end{Corollary}

Corollary \ref{cor:lem:gs-expansion} is the extension of Burchnall's formula \cite[equation~(5)]{Burc} for the Hermite polynomials, see~\eqref{eq:exCOP-Hermite-Burchnall2}.

We can next use Theorem \ref{lem:gs-expansion} in \eqref{eq:gs-LRadjointiterated} to get
\begin{gather*}
\sum_{k=0}^n \al^n_k \int_\R \frac{\eta^k (w_{\nu+k\si})(x)}{w_{\nu}(x)}\eta^k \big(p_{n-k}^{(\nu+k\si)}\big)(x)(T_{k,n}f)(x) g(x) d\mu^{(\nu)}(x) \\
 \qquad{} =\int_\R f(x) \bigl( L_{\nu+(n-1)\si} L_{\nu+(n-2)\si} \cdots L_{\nu+\si} L_{\nu} g\bigr)(x) d\mu^{(\nu+n\si)}(x).
\end{gather*}
Here we assume that each of the integrals converges. In particular, assuming that $L_\nu =\cD$ is a~lowering operator independent of~$\nu$, we obtain Corollary~\ref{cor:lem:gs-expansion2}.

\begin{Corollary}\label{cor:lem:gs-expansion2} Assuming the convergence of the integrals and with the notation and assumptions as above, we have
\begin{gather*}
\sum_{k=0}^n \al^n_k \int_\R \frac{\eta^k (w_{\nu+k\si})(x)}{w_{\nu}(x)} \eta^k \big(p_{n-k}^{(\nu+k\si)}\big)(x)(T_{k,n}f)(x) g(x) d\mu^{(\nu)}(x)\\
\qquad{} =\int_\R f(x) \bigl( \cD^n g\bigr)(x) d\mu^{(\nu+n\si)}(x).
\end{gather*}
\end{Corollary}

Corollary \ref{cor:lem:gs-expansion2} in case of the Hermite polynomials is not contained in Burchnall's paper~\cite{Burc}. Note that for $g(x)=x^p$ with $p<n$ the right hand side of Corollary~\ref{cor:lem:gs-expansion2} is zero.

\section{Example: classical orthogonal polynomials}\label{sec:examples-classicalOPS}

For the classical orthogonal polynomials in the Askey scheme, see~\cite{KoekLS,KoekS}, i.e., for the Jacobi, Laguerre and Hermite polynomials, we see that in Section~\ref{sec:general-setup} all the assumptions are fulfilled. Moreover, we can take $\cD= \partial = \frac{d}{dx}$, so that Leibniz formula is the usual one with $\al^n_k=\binom{n}{k}$, $\eta$~is the identity, and $T_{k,n} = \frac{d^k}{dx^k} = \cD^k$ is a power of the lowering operator.

\subsection{Example: Hermite polynomials}\label{ssec:Hermite}

In this case we discuss Burchnall's motivating example of the Hermite polynomials \cite{Burc}, and we show how to extend some of Burchnall's results using the generating function~\eqref{eq:Hermitegenerating} for Hermite polynomials. In this case the parameter set $\cV=\{0\}$ as in Section~\ref{sec:general-setup} is trivial. So $\si = 0$, and $w(x)=\exp\big({-}x^2\big)$ is independent of $\nu$, and so is $R_\nu = M_w^{-1}\circ \frac{d}{dx}\circ M_w = \big(\frac{d}{dx}-2x\big)$. Then the polynomials of Section~\ref{sec:general-setup} identify with the standard Hermite polynomials up to a sign; $(-1)^n H_n(x) = p_n^{(\nu)}(x)$, where $H_n(x) = (2x)^n {}_2F_0\big({-}n/2,-(n-1)/2;-;-x^{-1}\big)$. Theorem \ref{lem:gs-expansion} then gives \eqref{eq:exCOP-Hermite-Burchnall}, which is \cite[equation~(3)]{Burc}. Corollary~\ref{cor:lem:gs-expansion} then gives \eqref{eq:exCOP-Hermite-Burchnall2}, which was derived by Burchnall \cite[equation~(5)]{Burc} and Nielsen \cite[p.~22]{Niel}. For this derivation we use that $\frac{d}{dx}\big(\frac{d}{dx}-2x\big)x^{n-1} = -2nx^{n-1} + \text{l.o.t.}$
so that $\ell_n=-2n$, or $\frac{dH_n}{dx}(x) = 2n H_{n-1}(x)$.

For the Hermite polynomials, Burchnall's identity \eqref{eq:exCOP-Hermite-Burchnall} only involves the Hermite polyno\-mials. This gives the opportunity to elaborate a~bit more on Burchnall's identity. Multiply \eqref{eq:exCOP-Hermite-Burchnall} by $t^n/n!$ and sum over $n\in \N$ to obtain
\begin{gather}
\left( \exp\left(t \left(\frac{d}{dx} -2x \right) \right) f\right)(x) = \sum_{n=0}^\infty \sum_{k=0}^n \frac{t^k}{k!} f^{(k)}(x) \frac{(-t)^{n-k}}{(n-k)!} H_{n-k}(x) \nonumber\\
\hphantom{\left( \exp\left(t \left(\frac{d}{dx} -2x \right) \right) f\right)(x)}{} = f(x+t) \exp\big({-}2xt-t^2\big)\label{eqZ-B}
\end{gather}
after interchanging summation and using the generating function
\begin{gather}\label{eq:Hermitegenerating}
\sum_{n=0}^\infty H_n(x) \frac{t^n}{n!} = \exp\big(2xt-t^2\big)
\end{gather}
for the Hermite polynomials, see, e.g., \cite[equation~(6.1.7)]{AndrAR}, \cite[equation~(9.15.10)]{KoekLS}, \cite[equation~(1.13.10)]{KoekS}. Naturally we have to assume that~$f$ is sufficiently smooth, say real analytic, and~$t$ sufficiently small in~\eqref{eqZ-B}.

\begin{Remark} \label{rmk:Zassenhaus}
To see that \eqref{eqZ-B} follows from the Zassenhaus formula, or the inverse to the Baker--Campbell--Hausdorff formula, let $A = -2x$ (viewed as multiplication operator), $B = \frac{d}{dx}$, so that $[A, B] = 2$ and $[A,[A,B]]=0$, $[B,[A,B]]=0$, and hence all higher order commutators vanish as well. Zassenhaus's formula is, see \cite{CasaMN}, \cite[Section~4]{Magn},
\begin{gather*}
\exp(t(A+B)) = \exp(tA) \exp(tB) \exp\big(t^2 Z_2(A,B)\big) \exp\big(t^3 Z_3(A,B)\big) \cdots,
\end{gather*}
where $Z_m(A,B)$ is a homogeneous polynomial of degree $m$ in the non-commuting variables $A$ and $B$ given in terms of commutators;
\begin{gather*} 
Z_2(A, B) = - \frac{1}{2} [A, B], \qquad Z_3(A, B) = - \frac{1}{3} [[A, B],B] - \frac{1}{6} [[A, B],A], \qquad \dots.
\end{gather*}
So for this choice of $A$ and $B$ we have $Z_m(A,B)=0$ for $m\geq 3$ and $Z_2(A,B)= -1$. Hence, $\exp\bigl(t \big(\frac{d}{dx} -2x \big) \bigr)= \exp(-2xt) \exp\big(t\frac{d}{dx}\big) \exp\big({-}t^2\big)$, which gives the same result as~\eqref{eqZ-B} when acting on~$f$.
\end{Remark}

Note that Corollary \ref{cor:lem:gs-expansion2} after multiplication by $n!$ gives
\begin{gather*}
\int_\R e^{-x^2} \frac{g^{(n)} (x)}{n!} f(x) dx = \int_\R e^{-x^2} g(x) \sum_{k=0}^n \frac{(-1)^kf^{(k)} (x)}{k!} \frac{H_{n-k}(x)}{(n-k)!} dx.
\end{gather*}
Multiplying by $t^n$, summing over $n$ and using the generating function~\eqref{eq:Hermitegenerating} gives
\begin{gather}\label{eq:Hermite-intversion-shift}
\int_\R e^{-x^2} g(x+t) f(x) dx = e^{-t^2} \int_\R e^{2xt -x^2} f(x-t) g(x) dx,
\end{gather}
which follows directly by a change of variables. So, we can view \eqref{eq:Hermite-intversion-shift} as an integrated version of Burchnall's identity \eqref{eq:exCOP-Hermite-Burchnall}. From Corollary~\ref{cor:lem:gs-expansion2} with $f(x) =\exp(-xt)$ we see that $\sum_{k=0}^n \binom{n}{k} H_{n-k}(x) (-t)^k$ form orthogonal polynomials with respect to $e^{-xt} e^{-x^2}$. On the other hand, from the integrated version~\eqref{eq:Hermite-intversion-shift} we see that $H_n\big(x-\frac12 t\big)$ also form orthogonal polynomials with respect to $e^{-xt} e^{-x^2}$. So they are equal, up to a constant which can be determined by considering leading coefficients. This gives
\begin{gather}\label{eq:Hermite-Toda-OPs}
H_n(x-\frac12 t) = \sum_{k=0}^n (-t)^k \binom{n}{k} H_{n-k}(x)
\end{gather}
as the orthogonal polynomials with respect to $e^{-x^2} e^{-xt}$. Note that \eqref{eq:Hermite-Toda-OPs} can directly be established using the generating function~\eqref{eq:Hermitegenerating}. This gives the well-known Hermite case of Proposition~\ref{prop:explicit-Toda-solutions}. Note that~\eqref{eq:Hermite-Toda-OPs} is a special case of a more general convolution identity \cite[Corollary~3.8]{KoelVdJ}.

\subsection{Example: Laguerre polynomials}\label{ssec:Laguerre}

In this case, in the notation of Section~\ref{sec:general-setup}, $\cV= \{ \nu\in \R\,|\, \nu>-1\}\subset \R$, $\si = 1\in \R$. The measure $w_\nu(x)= x^\nu e^{-x}$ is absolutely continuous with respect to the Lebesgue measure on $[0,\infty)$. The polynomials $p^{(\nu)}_n(x)$ correspond to $n! L^{(\nu)}_n(x)$ by comparing with \cite[equation~(1.11.9)]{KoekS}. Here the Laguerre polynomials are defined by
\begin{gather*}
L^{(\nu)}_n(x) = \frac{(\nu+1)_n}{n!} \rFs{1}{1}{-n}{\nu+1}{x}.
\end{gather*}

In this case $R_\nu = M_\nu^{-1}\circ \frac{d}{dx}\circ M_{\nu+1} = \nu+1 + x\frac{d}{dx}$
and by Theorem \ref{lem:gs-expansion} we have
\begin{gather}\label{eq:Laguerre-lemmagsexpansion}
(R_\nu R_{\nu+1} \cdots R_{\nu+n-1} f)(x) = \sum_{k=0}^n
\frac{n!}{k!} f^{(k)}(x) x^k L_{n-k}^{(\nu+k)}(x),
\end{gather}
see Carlitz \cite[equation~(4)]{Carl}, where it is proved by induction on~$n$. Corollary \ref{cor:lem:gs-expansion} gives \eqref{eq:Laguerre-corlemmagsexpansion}, which is \cite[equation~(7)]{Carl}. Carlitz then continues to prove from \eqref{eq:Laguerre-corlemmagsexpansion} an extension of the gene\-ra\-ting function \eqref{eq:Laguerregenerating}, which was previously obtained in Rainville \cite[Section~119, equation~(9)]{Rain}. Rainville's generating function is equivalent to~\eqref{eq:Laguerre-corlemmagsexpansion}.

Corollary \ref{cor:lem:gs-expansion2} yields
\begin{gather}\label{eq:Laguerre-cor2lemmagsexpansion}
\int_0^\infty x^{\nu+n} e^{-x} g^{(n)}(x) f(x) dx = (-1)^n \int_0^\infty x^\nu e^{-x} g(x)\sum_{k=0}^n\frac{n!}{k!} f^{(k)}(x) x^{k} L_{n-k}^{(\nu+k)}(x) dx
\end{gather}
for $\nu>-1$. Multiplying \eqref{eq:Laguerre-cor2lemmagsexpansion} by~$u^n$, dividing by $(-1)^n n!$, summing over $n$ and using the generating function
\begin{gather}\label{eq:Laguerregenerating}
\sum_{n=0}^\infty u^n L_n^{(\nu)}(x) = (1-u)^{-\nu-1}\exp(xu/(u-1)), \qquad |u|<1,
\end{gather}
see, e.g., \cite[equation~(1.11.10)]{KoekS}, gives
\begin{gather}\label{eq:Laguerre-intversion-shift}
(1-u)^{\nu+1} \int_0^\infty x^\nu e^{-x} f(x) g(x-ux) dx =\int_0^\infty x^\nu e^{-x} \exp\left(\frac{-xu}{1-u}\right) g(x) f\left(\frac{x}{1-u}\right) dx\!\!\!
\end{gather}
for sufficiently smooth functions $f$ and $g$ and $-1<u<1$. Note that \eqref{eq:Laguerre-intversion-shift} is the integrated version, and it can be proved directly.

From \eqref{eq:Laguerre-intversion-shift}, i.e., Corollary \ref{cor:lem:gs-expansion2} for the Laguerre polynomials, we can obtain an expression for the orthogonal polynomials with respect to $e^{-xt} x^\nu e^{-x}$ on $[0,\infty)$ by taking $f(x) =\exp(-xt)$ in~\eqref{eq:Laguerre-cor2lemmagsexpansion}, and take the polynomial $g(x) = x^p$, with $p<n$. This gives the right hand side of
\begin{gather}\label{eq:Laguerre-Toda-OPs}
L_n^{(\nu)}(x(1+t)) = \sum_{k=0}^n \frac{1}{k!} (-t)^k x^{k} L_{n-k}^{(\nu+k)}(x)
\end{gather}
as the orthogonal polynomials with respect to $e^{-xt} x^\nu e^{-x}$ on $[0,\infty)$. On the other hand, by a~straightforward calculation using~\eqref{eq:Laguerre-intversion-shift}, the polynomials on the left hand side are also orthogonal with respect to the same weight. Hence, \eqref{eq:Laguerre-Toda-OPs} follows up to a constant. This constant is determined by evaluating at $x=0$. The result~\eqref{eq:Laguerre-Toda-OPs} is a convolution type identity, and can be directly proved using the generating function~\eqref{eq:Laguerregenerating}. See \cite[Section~3]{KoelVdJ} for generalizations of convolution identities for Laguerre polynomials.

In particular, we find the solutions to the Toda equation~\eqref{eq:intro-Toda} for the Laguerre case, see the Laguerre case of Proposition~\ref{prop:explicit-Toda-solutions}.

\begin{Remark} Gould and Hopper \cite{GoulH} have generalized the Hermite and Laguerre cases simultaneously. In particular, Gould and Hopper derive an operational formula for the associated differential operator \cite[Section~4]{GoulH}. Al-Salam \cite[p.~130]{AlSa} shows that the operational identities of Carlitz \cite{Carl} and Gould--Hopper \cite{GoulH} are equivalent.
\end{Remark}

\subsection{Example: Jacobi polynomials}\label{ssec:Jacobi}

In this case the ingredients of Section~\ref{sec:general-setup} correspond to $\cV= \{ (\al,\be)\in \R^2 \,|\, \al>-1, \be>-1\}\subset \R$, $\si = (1,1)\in \R^2$. The measure is $w_{\al,\be}(x) = (1-x)^\al (1+x)^\be$ with respect to Lebesgue measure on $[-1,1]$. The polynomials $p^{(\nu)}_n(x)$ of Section~\ref{sec:general-setup} correspond to $(-2)^n n! P^{(\al,\be)}_n(x)$ by comparing with \cite[equation~(1.8.9)]{KoekS}. Here the Jacobi polynomials are defined by
\begin{gather*}
P^{(\al,\be)}_n(x) = \frac{(\al+1)_n}{n!} \rFs{2}{1}{-n, n+\al+\be+1}{\al+1}{\frac12(1-x)}.
\end{gather*}
Then
\begin{gather*}
R_{\al,\be} = M_{\al,\be}^{-1}\circ \frac{d}{dx} \circ M_{\al+1,\be+1} = \big(1-x^2\big)\frac{d}{dx} + \bigl( \be-\al - x(\al+\be+2)\bigr),
\end{gather*}
and Theorem \ref{lem:gs-expansion} gives
\begin{gather}
(R_{\al, \be}R_{\al+1, \be+1}\cdots R_{\al+n-1, \be+n-1} f)(x)\nonumber\\
\qquad{} = \sum_{k=0}^n \frac{n!}{k!} (-2)^{n-k} f^{(k)}(x) \big(1-x^2\big)^k P_{n-k}^{(\al+k, \be+k)}(x)\label{eq:Jacobi-lemmagsexpansion}
\end{gather}
for $f\in C^\infty(\R)$. Corollary~\ref{cor:lem:gs-expansion} gives
\begin{gather}
 \binom{n+m}{n} P_{m+n}^{(\al, \be)}(x) = \sum_{k=0}^n \frac{1}{k!} \left( \frac{1-x^2}{-4}\right)^k(\al+\be+2n+m+1)_k\nonumber\\
 \hphantom{\binom{n+m}{n} P_{m+n}^{(\al, \be)}(x) =}{} \times P_{n-k}^{(\al+k, \be+k)}(x) P_{m-k}^{(\al+n+k, \be+n+k)}(x),\label{eq:Jacobi-corlemmagsexpansion}
\end{gather}
since $\frac{dP_n^{(\al,\be)}}{dx}(x) = \frac12 (n+\al+\be+1) P_{n-1}^{(\al+1,\be+1)}(x)$, see, e.g., \cite[equation~(1.8.6)]{KoekS}. Singh~\cite{Sing} obtains the operational formula~\eqref{eq:Jacobi-lemmagsexpansion} and expansion~\eqref{eq:Jacobi-corlemmagsexpansion}.

The integrated version follows from Corollary \ref{cor:lem:gs-expansion2}:
\begin{gather*}
\int_{-1}^1 (1-x)^{\al+n} (1+x)^{\be+n} g^{(n)}(x) f(x) dx \nonumber \\
\qquad{} = \int_{-1}^1 (1-x)^{\al} (1+x)^{\be} g(x) \sum_{k=0}^n\frac{n!}{k!}f^{(k)}(x)(-2)^{n-k} \big(1-x^2\big)^k P_{n-k}^{(\al+k, \be+k)}(x) dx.
\end{gather*}
By taking $f(x) = \exp(-xt)$, $g(x)=x^p$ we see that
\begin{gather*}
\int_{-1}^1 e^{-xt} (1-x)^{\al} (1+x)^{\be} x^p\\
\qquad{}\times \sum_{k=0}^n\frac{n!}{k!} (-t)^k(-2)^{n-k} (n-k)! \big(1-x^2\big)^k P_{n-k}^{(\al+k, \be+k)}(x) dx = 0, \qquad p < n,
\end{gather*}
so that we find only partial orthogonality, since the sum over $k$ gives a polynomial of deg\-ree~$2n$, instead of $n$. For a detailed study of the orthogonal polynomials with respect to the Toda modification of the Jacobi weight $e^{-xt} (1-x)^{\al} (1+x)^{\be}$ on $[-1,1]$, we refer to Basor, Chen and Ehrhardt~\cite{BasoCE}, where the relation to the Painlev\'e V equation is discussed.

\section[The backward shift operator $\nabla$]{The backward shift operator $\boldsymbol{\nabla}$}\label{sec:backwardshift}

In this section we consider the results of Section~\ref{sec:general-setup} in case the derivative $\partial$ is given by the backward shift operator $\nabla$, defined by $\bigl( \nabla f\bigr) (x) = f(x)-f(x-1)$. In the Askey scheme, this corresponds to the families of Meixner and Charlier polynomials. The Krawtchouk polynomials form a family of finite discrete orthogonal polynomials, also contained in the part of the Askey scheme corresponding to $\nabla$. We leave the Krawtchouk case to the reader, but we include the case in Proposition~\ref{prop:explicit-Toda-solutions} for completeness. In order to apply the results of Section~\ref{sec:general-setup} we need to have the Leibniz formula \eqref{eq:gs-Leibniz} explicitly:
\begin{gather}\label{eq:backwardshift-Leibniz}
\nabla^n (fg)(x) = \sum_{k=0}^n \binom{n}{k} \big(\nabla^kf\big)(x)\big(S^k \nabla^{n-k}g\big)(x) = \sum_{k=0}^n \binom{n}{k} \big(\nabla^{n-k}f\big)(x)\big(S^{n-k} \nabla^{k}g\big)(x),\!\!\!\!
\end{gather}
where $Sf(x) = f(x-1)$. This follows by induction on $n$, since $S$ and $\nabla$ commute. Note that, upon comparing~\eqref{eq:backwardshift-Leibniz} with~\eqref{eq:gs-Leibniz}, we see that there are two choices for the homomorphism $\eta$ of Section~\ref{sec:general-setup}; $\eta=1$ or $\eta=S$. Different choices lead to different expansions.

\subsection{The Meixner polynomials}\label{sssec:Toda-backwardshift-Meixner}

In the notation of Section~\ref{sec:general-setup} we have $\cV =\{ (\be, c) \in \R^2 \,|\, \be >0, 0<c<1\}$, and we let $\si=(1,0)$, then $\nu\in \cV$ implies $\nu+n\si\in \cV$ for all $n\in \N$. Put
\begin{gather*} 
w(x;\be,c) = \frac{(\be)_x}{x!}c^x, \qquad x\in \N,
\end{gather*}
so that the orthogonality measure is a discrete measure supported in $\N$ with weight at $x\in \N$ equal to $w(x;\be, c)$. Then the raising operator $R_\nu$ corresponds to $R_{\be,c}= M_{\be,c}^{-1}\circ \nabla \circ M_{\be+1,c}$, so
\begin{gather*}
\bigl( R_{\be,c} f\bigr) (x) = \frac{\be+x}{\be} f(x) - \frac{x}{c \be} f(x-1).
\end{gather*}
By \cite[equation~(1.9.9)]{KoekS} $p^{(\nu)}_n(x)$ corresponds precisely to the Meixner polynomials $M_n(x;\be,c)$ defined by
\begin{gather*}
M_n(x;\be,c) = \rFs{2}{1}{-n, -x}{\be}{\frac{c-1}{c}}.
\end{gather*}
The adjoint $L_\nu$ corresponds to $-\Delta$ is independent of $\nu$, where $\Delta f(x) = f(x+1)-f(x)$ is the forward shift operator.

\subsubsection[The Leibniz rule with $\eta=1$]{The Leibniz rule with $\boldsymbol{\eta=1}$}\label{ssssec:backwardshift-Meixner-etais1}

Comparing \eqref{eq:backwardshift-Leibniz} with \eqref{eq:gs-Leibniz} we take $\partial=\nabla$, $\eta=1$, $\al^n_k=\binom{n}{k}$, and $T_{k,n} = S^{n-k} \nabla^{k} = (-1)^k S^n \Delta^k$. Then we identify $\frac{\eta^k (w_{\nu+k\si})(x)}{w_{\nu}(x)}$ with
\begin{gather*}
\frac{w(x; \be+k,c)}{w(x;\be,c)} = \frac{(\be+k)_x}{(\be)_x} = \frac{(\be+x)_k}{(\be)_k}.
\end{gather*}
The operational form of Theorem \ref{lem:gs-expansion} then gives
\begin{gather}\label{eq:backwardshift-etais1-operational}
\bigl( R_{\be,c} R_{\be,c}\cdots R_{\be+n-1,c}f\bigr)(x) = \sum_{k=0}^n \binom{n}{k} \frac{(\be+x)_k}{(\be)_k} (-1)^k M_{n-k}(x;\be+k,c) \bigl( S^n\Delta^k f\bigr)(x).\!\!\!\!
\end{gather}

By \cite[equation~(1.9.7)]{KoekS}, $T_{k,n} p^{(\nu+n\si)}_m(x)$ corresponds to
\begin{gather*}
(-1)^k \frac{(m-k+1)_k}{(\be+n)_k} \left( \frac{c-1}{c}\right)^k M_{m-k}(x-n; \be+n+k,c).
\end{gather*}
Corollary \ref{cor:lem:gs-expansion} gives, after a simplification,
\begin{gather}
M_{m+n}(x;\be,c) = \sum_{k=0}^{n \wedge m}\frac{(-n)_k(-m)_k (\be+x)_k}{k! (\be)_k (\be+n)_k}\left( \frac{c-1}{c}\right)^k \nonumber\\
\hphantom{M_{m+n}(x;\be,c) =}{} \times M_{n-k}(x;\be+k,c) M_{m-k}(x-n; \be+n+k,c).\label{eq:Toda-backwardshift-Meixner-expansion-etais1}
\end{gather}

Applying Corollary \ref{cor:lem:gs-expansion2} with $f(x) = e^{-xt}$, which is a joint eigenfunction of~$\nabla$, $\Delta$ and $S$, and $g(x)=x^p$, $p<n$, we see that the polynomials of degree~$n$
\begin{gather*}
\sum_{k=0}^n \binom{n}{k} \frac{(\be+x)_k}{(\be)_k} M_{n-k}(x;\be+k,c) e^{nt}\big(1-e^{-t}\big)^k
\end{gather*}
are orthogonal with respect to the discrete measure $h\mapsto \sum_{x=0}^ \infty w(x;\be,c) e^{-xt} h(x)$, for which all moments exist for $t> \ln(c)$. Since this measure is again the measure for the Meixner polynomials with
parameters $(\be, ce^{-t})$ we find
\begin{gather}\label{eq:Toda-Meixner-expansionorthpols-etais1}
M_n\big(x;\be, ce^{-t}\big) =\sum_{k=0}^n \frac{(-n)_k (\be+x)_k }{k!(\be)_k } M_{n-k}(x;\be+k, c) e^{nt}\big(1-e^{-t}\big)^k,
\end{gather}
where the constant is obtained by evaluating at $x=0$, since $M_n(0;\be,c)=1$. Note that~\eqref{eq:Toda-Meixner-expansionorthpols-etais1} is a convolution identity, which can be obtained directly from the generating function, see, e.g., \cite[equation~(1.9.11)]{KoekS},
\begin{gather}\label{eq:Meixner-generating}
\sum_{n=0}^\infty \frac{(\be)_n}{n!} t^n M_n(x;\be,c) =\left(1-\frac{t}{c}\right)^x (1-t)^{-x-\be}
\end{gather}
and the binomial theorem. More general convolution identities for Meixner polynomials can be found in \cite[Section~3]{KoelVdJ}.

In particular, we find an explicit solution to the Toda lattice \eqref{eq:intro-Toda} by looking at the three-term recurrence relation for the monic version of $M_n(c;\be, ce^{-t})$. This in particular gives the explicit solution for the Meixner case of Proposition~\ref{prop:explicit-Toda-solutions}.

\subsubsection[The Leibniz rule with $\eta=S$]{The Leibniz rule with $\boldsymbol{\eta=S}$}\label{ssssec:backwardshift-Meixner-etaisS}
Comparing \eqref{eq:backwardshift-Leibniz} with \eqref{eq:gs-Leibniz} we can also take $\partial=\nabla$, $\eta=S$, $\al^n_k=\binom{n}{k}$, and $T_{k,n} = \nabla^{k}$ in the general set-up. Then we identify $\frac{\eta^k (w_{\nu+k\si})(x)}{w_{\nu}(x)}$ with
\begin{gather*}
\frac{w(x-k; \be+k,c)}{w(x;\be,c)} = \frac{(x-k+1)_k}{(\be)_k c^k} = \frac{(-1)^k (-x)_k}{(\be)_k c^k}.
\end{gather*}
The operational form of Theorem~\ref{lem:gs-expansion} is{\samepage
\begin{gather*}
\bigl( R_{\be,c} R_{\be,c}\cdots R_{\be+n-1,c}f\bigr)(x) = \sum_{k=0}^n \binom{n}{k} \frac{(-x)_k}{(\be)_k} (-c)^{-k}M_{n-k}(x-k;\be+k,c) \bigl( \nabla^k f\bigr)(x),
\end{gather*}
which gives an alternative expansion for the same operator as in \eqref{eq:backwardshift-etais1-operational}.}

Rewrite $T_{k,n} = \nabla^k = (-1)^k S^k \Delta^k$, so that by \cite[equation~(1.9.7)]{KoekS} the polynomial $T_{k,n} p_m^{(\nu+n\si)}$ is identified with $(-1)^k \frac{(m-k+1)_k}{(\be+n)_k} \big( \frac{c-1}{c}\big)^k M_{m-k}(x-k; \be+n+k,c)$. Corollary \ref{cor:lem:gs-expansion} yields, after simplification,
\begin{gather}
M_{m+n}(x;\be,c) = \sum_{k=0}^{n\wedge m} \frac{(-n)_k (-m)_k (-x)_k}{k! (\be)_k (\be+n)_k} \left( \frac{c-1}{c^2}\right)^k\nonumber\\
\hphantom{M_{m+n}(x;\be,c) =}{} \times M_{n-k}(x-k;\be+k,c) M_{m-k}(x-k; \be+n+k,c).\label{eq:Toda-backwardshift-Meixner-expansion-etaisS}
\end{gather}
Note that \eqref{eq:Toda-backwardshift-Meixner-expansion-etais1} and~\eqref{eq:Toda-backwardshift-Meixner-expansion-etaisS} give different expansions for the same Meixner polynomial.

Using Corollary \ref{cor:lem:gs-expansion2} with $f(x) = e^{-xt}$ we find that the orthogonal polynomials for the Toda modification have an expansion as in the right hand side. Since we know the corresponding orthogonal polynomials as Meixner poynomials with parameters~$(\be, ce^{-t})$ we get
\begin{gather}\label{eq:Toda-Meixner-expansionorthpols-etaiseta}
M_n(x;\be, ce^{-t}) = \sum_{k=0}^n \frac{(-n)_k (-x)_k}{k! (\be)_k c^k}M_{n-k}(x-k;\be+k,c) \big(1-e^t\big)^k,
\end{gather}
where the constant follows by comparing leading coefficients using \cite[equation~(1.9.4)]{KoekS} and the binomial theorem. Note that \eqref{eq:Toda-Meixner-expansionorthpols-etais1} and~\eqref{eq:Toda-Meixner-expansionorthpols-etaiseta} give different expansions for the same Meixner polynomial. Again, the convolution identity \eqref{eq:Toda-Meixner-expansionorthpols-etaiseta} can be proved directly using the generating function \eqref{eq:Meixner-generating} and the binomial formula, and see \cite[Corollary~3.6]{KoelVdJ} for generalizations.

\subsection{The Charlier polynomials}\label{sssec:Toda-backwardshift-Charlier}

The orthogonality measure for the Charlier polynomials is a discrete measure supported on $\N$ and with corresponding weights $w_a(x) = \frac{a^x}{x!}$, $x\in \N$. In the notation of Section~\ref{sec:general-setup} we have $\cV=\{ a>0\} \subset \R$. In this case the raising operator is independent of $a$, so that $\si=0$. So $R_a = M_a^{-1}\circ \nabla \circ M_a$, so
\begin{gather*}
\bigl( R_a f\bigr) (x) = f(x) - \frac{x}{a} f(x-1)
\end{gather*}
and $L = \cD = -\De$. Then the polynomials $p_n^{(\nu)}$ correspond precisely to the Charlier polynomials $C_n(x;a) = {}_2F_0\big({-}n, -x;-;\frac{-1}{a}\big)$ as in \cite[Section~1.12]{KoekS}. It should be noted that these results can also be obtained by a limit transition from the corresponding results for the Meixner polynomials.

\subsubsection[The Leibniz rule with $\eta=1$]{The Leibniz rule with $\boldsymbol{\eta=1}$}\label{ssssec:Toda-backwardshift-Charlier-etais1}

Here we give the analogues as in Section~\ref{ssssec:backwardshift-Meixner-etais1}. So take $\partial=\nabla$, $\eta=1$, $\al^n_k=\binom{n}{k}$, and $T_{k,n} = S^{n-k} \nabla^{k} = (-1)^k S^n \Delta^k$. Then we identify $\frac{\eta^k (w_{\nu+k\si})(x)}{w_{\nu}(x)}$ with~$1$. So the operational formula of Theorem~\ref{lem:gs-expansion} gives
\begin{gather}\label{eq:Charlier-operational-expansion-etais1}
\bigl( R_a^n f\bigr)(x) = \sum_{k=0}^n \binom{n}{k} C_{n-k}(x; a) \bigl( S^{n-k} \nabla^k f\bigr)(x).
\end{gather}
By \cite[equation~(1.9.7)]{KoekS} $T_{k,n} p^{(\nu+n\si)}_m(x)$ corresponds to $(-1)^k \frac{(-m)_k}{a^k} C_{m-k}(x-n; a)$, and Corollary~\ref{cor:lem:gs-expansion} gives
\begin{gather}\label{eq:Toda-backwardshift-Charlier-expansion-etais1}
C_{n+m}(x;a) = \sum_{k=0}^{n \wedge m} \frac{(-n)_k (-m)_k}{k! a^k} C_{n-k}(x;a) C_{m-k}(x-n;a).
\end{gather}

Using Corollary \ref{cor:lem:gs-expansion2} with $f(x) = e^{-xt}$ we find an explicit expansion for the polynomials orthogonal with respect to $w_a(x) e^{-xt}$ for $x\in \N$. Since these are the polynomials $C_n(x;ae^{-t})$ we obtain, cf.\ Section~\ref{ssssec:backwardshift-Meixner-etais1},
\begin{gather}\label{eq:Toda-Charlier-expansionorthpols-etais1}
\sum_{k=0}^n \frac{(-n)_k}{k!} C_{n-k}(x;a) e^{nt} \big(1-e^{-t}\big)^k = C_n\big(x;ae^{-t}\big).
\end{gather}
This is a convolution type identity, which can be obtained from the generating function for Charlier polynomials; $\sum_{n=0}^\infty \frac{z^n}{n!} C_n(x;a) = e^z \big(1-\frac{z}{a}\big)^x$, see, e.g., \cite[equation~(1.12.11)]{KoekS}, and the Taylor expansion for the exponential function. It also follows from \eqref{eq:Toda-Meixner-expansionorthpols-etais1} by taking the limit from Meixner to Charlier polynomials, see \cite[equation~(2.9.2)]{KoekS}.

In particular, we find an explicit solution to the Toda lattice \eqref{eq:intro-Toda} by looking at the three-term recurrence relation for the monic version of~$C_n(x;ae^{-t})$. This in particular gives the explicit solution for the Charlier case of Proposition~\ref{prop:explicit-Toda-solutions}.

\subsubsection[The Leibniz rule with $\eta=S$]{The Leibniz rule with $\boldsymbol{\eta=S}$}\label{ssssec:Toda-backwardshift-Charlier-etaisS}
As in Section~\ref{ssssec:backwardshift-Meixner-etaisS} we can also take $\partial=\nabla$, $\eta=S$, $\al^n_k=\binom{n}{k}$, and $T_{k,n} = \nabla^{k}$ in the general set-up. Then we identify $\frac{\eta^k (w_{\nu+k\si})(x)}{w_{\nu}(x)}$ with $(-a)^{-k} (-x)_k$. So the operational formula of Theo\-rem~\ref{lem:gs-expansion} gives the alternative
\begin{gather*}
\bigl( R_a^n f\bigr)(x) = \sum_{k=0}^n \binom{n}{k}(-a)^{-k} (-x)_k C_{n-k}(x-k; a) \bigl( \nabla^k f\bigr)(x)
\end{gather*}
to \eqref{eq:Charlier-operational-expansion-etais1}.

Then $T_{k,n} = \nabla^k = (-1)^k S^k \Delta^k$, so that by \cite[equation~(1.9.7)]{KoekS} the polynomial $T_{k,n} p_m^{(\nu+n\si)}$ is identified with $(-a)^{-k} (-m)_k C_{m-k}(x-k; a)$. So Corollary~\ref{cor:lem:gs-expansion} gives
\begin{gather}\label{eq:Toda-backwardshift-Charlier-expansion-etaisS}
C_{n+m}(x;a) = \sum_{k=0}^{n \wedge m} \frac{(-n)_k (-m)_k}{k!} \big({-}a^2\big)^{-k} C_{n-k}(x-k;a) C_{m-k}(x-k;a).
\end{gather}
Note that \eqref{eq:Toda-backwardshift-Charlier-expansion-etais1} and \eqref{eq:Toda-backwardshift-Charlier-expansion-etaisS} give different expansions.

Similarly, using Corollary \ref{cor:lem:gs-expansion2} with $f(x) = e^{-xt}$ we get in a similar way the equality
\begin{gather}\label{eq:Toda-Charlier-expansionorthpols-etaiseta}
\sum_{k=0}^n \frac{(-n)_k(-x)_k}{k!} C_{n-k}(x-k;a) a^{-k} \big(1-e^{t}\big)^k = C_n\big(x;ae^{-t}\big)
\end{gather}
up to a constant which can be determined by considering leading coefficients. Then \eqref{eq:Toda-Charlier-expansionorthpols-etaiseta} can be proved using the same generating function for the Charlier polynomials and the binomial theorem. Again, \eqref{eq:Toda-Charlier-expansionorthpols-etais1} and~\eqref{eq:Toda-Charlier-expansionorthpols-etaiseta} give different exansions for the same polynomial.

\section[The difference operator $\frac{\de}{\de x}$]{The difference operator $\boldsymbol{\frac{\de}{\de x}}$}\label{sec:div-difference}

The difference operator $\frac{\de}{\de x}$ is defined by
\begin{gather*}
\frac{\de f}{\de x}(x) = \frac{1}{i} \left( f\left(x+\frac12 i\right)- f\left(x-\frac12 i\right)\right)
\end{gather*}
and the polynomials in the Askey scheme corresponding to this operator as part of the raising operator are the families of the Meixner--Pollaczek and the continuous Hahn polynomials. We discuss the Meixner--Pollaczek polynomials in detail.

In order to apply the results of Section~\ref{sec:general-setup} we need to have the Leibniz formula \eqref{eq:gs-Leibniz} explicitly:
\begin{gather}\label{eq:Toda-difference-Leibniz}
\left(\frac{\de}{\de x}\right)^n (fg)(x) = \sum_{k=0}^n \binom{n}{k} \left( \big(S^+\big)^{n-k} \left(\frac{\de}{\de x}\right)^{k}f\right)(x) \left( (S^-)^k \left(\frac{\de}{\de x}\right)^{n-k}g\right)(x),
\end{gather}
where $S^\pm f(x) = f\big(x\pm \frac12 i\big)$, assuming that the functions $f$ and $g$ are defined in a~sufficiently large strip $|\Im z|< \frac12 n+\ep$ for $\ep>0$. This follows by induction on $n$, since $S^\pm$ and~$\frac{\de}{\de x}$ commute. In this case the Leibniz formula is symmetric in shifting~$x$ to $x\pm\frac12 i$, so it suffices to consider only one of the two possibilities. We take $\eta= S^-$, $\al^n_k =\binom{n}{k}$, $\partial = \frac{\de}{\de x}$, $T_{k,n}= (S^+)^{n-k} \bigl( \frac{\de}{\de x}\bigr)^k$.

\subsection{The Meixner--Pollaczek polynomials}\label{sssec:Toda-difference-MP}

In this case $\cV =\{(\la,\vp) \,|\, \la>0, 0<\vp<\pi\}$, and $\si=\big(\frac12, 0\big)$ in the notation of Section~\ref{sec:general-setup}. The weight for the Meixner--Pollaczek polynomials
\begin{gather*}
w(x;\la, \vp) = \Ga(\la+ix) \Ga(\la-ix) e^{(2\vp-\pi)x}
\end{gather*}
is supported on $\R$. So $R_\nu$ corresponds to $R_{\la,\vp}= M^{-1}_{\la,\vp}\circ \frac{\de}{\de x} \circ M_{\la+\frac12,\vp}$, or
\begin{gather*}
\bigl(R_{\la,\vp} f\bigr)(x)= - e^{i\vp} (\la-ix) f\left(x+\frac12 i\right) - e^{-i\vp} (\la+ix) f\left(x-\frac12 i\right).
\end{gather*}

Comparing with \cite[equation~(1.7.9)]{KoekS} we find that $p^{(\nu)}_n$ corresponds to $(-1)^n n! P_n^{(\la)}(x;\vp)$, where the Meixner--Pollaczek polynomials are defined by
\begin{gather*}
P_n^{(\la)}(x;\vp) = \frac{(2\la)_n}{n!} e^{in\vp} \rFs{2}{1}{-n, \la+ix}{2\la}{1-e^{-2i\vp}}.
\end{gather*}
Note that $P_n^{(\la)}(x;\vp)= (-1)^n P_n^{(\la)}(-x;\pi-\vp)$ by Pfaff's transformation, see, e.g., \cite[Theorem~2.2.5]{AndrAR}, and this gives another explanation why we can restrict to one of the possibilities for the Leibniz formula~\eqref{eq:Toda-difference-Leibniz}. We identify $\frac{\eta^k (w_{\nu+k\si})(x)}{w_{\nu}(x)}$ with
\begin{gather*}
\frac{w\big(x-\frac{k}{2}i; \la +\frac12 k,\vp \big)}{w(x;\la,\vp )} = i^k e^{-ik\vp} (\la+ix)_k
\end{gather*}
so that the operational formula of Theorem \ref{lem:gs-expansion} gives
\begin{gather*}
\bigl(R_{\la,\phi} R_{\la+\frac12,\vp} \cdots R_{\la+\frac12 (n-1),\vp}f\bigr)(x)\\
\qquad{} = \sum_{k=0}^n \frac{n!}{k!} i^k e^{-ik\vp} (-1)^{n-k} (\la+ix)_k
P_{n-k}^{(\la+\frac{k}{2})}\left(x-\frac12k i;\vp\right) \left(\bigl(S^+\bigr)^{n-k} \left( \frac{\de}{\de x}\right)^k f\right)(x).\nonumber
\end{gather*}

Now $T_{k,n} p^{(\nu+n\si)}$ corresponds to $(-1)^m m! (2\sin\vp)^k P^{(\la+\frac12(n+k))}_{m-k}(x+\frac{i}{2}(n-k);\vp)$ by \cite[equation~(1.7.7)]{KoekS}, so that Corollary~\ref{cor:lem:gs-expansion} then gives, after simplification,
\begin{gather*}
\binom{n+m}{n} P^{(\la)}_{m+n}(x;\vp) =\sum_{k=0}^{n\wedge m}\frac{(-i)^k e^{-ik\vp}}{k!} (2\sin\vp)^k(\la+ix)_k \nonumber\\
\hphantom{\binom{n+m}{n} P^{(\la)}_{m+n}(x;\vp) =}{} \times P^{(\la+\frac12k)}_{n-k}\left(x-\frac{k}{2}i;\vp\right)
P^{(\la+\frac12(n+k))}_{m-k}\left(x+\frac{i}{2}(n-k);\vp\right).
\end{gather*}

Applying Corollary \ref{cor:lem:gs-expansion2} to $f(x)=e^{-xt}$ and $g(x) = x^p$, $0\leq p <n$, we find an expression for the orthogonal polynomials with respect to the measure $w(x;\la, \vp) e^{-xt}$ on $\R$. Since these polynomials are Meixner--Pollaczek polynomials for the weight $w\big(\cdot; \la, \vp-\frac12t\big)$, $2\vp-2\pi<t<2\vp$, we find the polynomials
\begin{gather}
P_n^{(\la)}\left(x;\vp-\frac12 t\right) \nonumber\\
\qquad{} = \sum_{k=0}^n \frac{i^k e^{-ik\vp}}{k!} (\la+ix)_k P^{(\la+\frac12 k)}_{n-k}\left(x-\frac{k}{2}i;\vp\right)\left(2\sin\left(\frac12 t\right)\right)^k e^{-\frac12 t(n-k)i},\label{eq:MP-Toda-expansion}
\end{gather}
where the constant follows by comparing the value at $x=i\la$. Note that~\eqref{eq:MP-Toda-expansion} is a convolution identity which can be directly obtained from the generating function for the Meixner--Pollaczek polynomials; $\sum_{n=0}^\infty u^n P^{(\la)}_n(x;\vp) = (1-e^{i\vp}u)^{-\la+ix}(1-e^{-i\vp}u)^{-\la-ix}$, see, e.g., \cite[equation~(1.7.11)]{KoekS}. For more general convolution identities for Meixner--Pollaczek polynomials, see \cite[Theorem~3.4]{KoelVdJ}. Having the orthogonal polynomials for the modified weight $w(x;\la, \vp) e^{-xt}$ on $\R$, the Toda equations \eqref{eq:intro-Toda} can be solved explicitly, see Proposition~\ref{prop:explicit-Toda-solutions}.

\section{Burchnall's identity for the Wilson polynomials}\label{sec:examples-Wilson}

As stated in Section~\ref{sec:general-setup} the Burchnall results can be derived for any polynomial family in the Askey scheme and its $q$-analogue. In this section we present the result for the top level of the Askey-scheme.

In this case $\cV=\{ \mathbf{a}=(a,b,c,d)\,|\, a>0, b> 0, c>0, d>0 \}\subset \R^4$ is the parameter space. The parameter space can be made more general, see, e.g., \cite[Section~1.1]{KoekS}, but the results can be directly extended to the more general parameter space by analytic extension in the parameter, and we stick to the real parameters. The shift corresponds to $\si=\big(\frac12,\frac12,\frac12,\frac12\big)$, see, e.g., \cite[equation~(1.1.7)]{KoekS}. The weight function is, up to a constant independent of $\mathbf{a}$ and a change of variable,
\begin{gather*}
\om (x; \mathbf{a}) = \om(x; a,b,c,d) \\
\hphantom{\om (x; \mathbf{a})}{} = \frac{\Ga(a+ix)\Ga(a-ix)\Ga(b+ix)\Ga(b-ix)\Ga(c+ix)\Ga(c-ix)\Ga(d+ix)\Ga(d-ix)}{2ix \Ga(2ix) \Ga(-2ix)}.
\end{gather*}
Then the Wilson polynomials $W_n(\cdot; \mathbf{a})$ satisfy the orthogonality
\begin{gather*}
 \int_0^\infty W_n\big(x^2;\mathbf{a}\big)W_m\big(x^2;\mathbf{a}\big) 2ix \om (x; \mathbf{a}) dx = 0, \qquad n\not=m, \\
 W_n\big(x^2;\mathbf{a}\big) = (a+b)_n (a+c)_n(a+d)_n \rFs{4}{3}{-n, n+a+b+c+d-1, a+ix, a-ix}{a+b, a+c, a+d}{1}.
\end{gather*}
So the raising operator $R_\nu$ in Section~\ref{sec:general-setup} corresponds to $R_{\mathbf{a}}=M^{-1}_{\mathbf{a}}\circ \frac{\de}{\de x^2} \circ M_{\mathbf{a}+\frac12}$, where $\mathbf{a}+\frac12= \big(a+\frac12,b+\frac12,c+\frac12,d+\frac12\big)$ and $\frac{\de}{\de x^2}f = \frac{1}{2ix}\bigl( f\big(x+\frac12 i\big)- f\big(x-\frac12 i\big)\bigr)$. Explicitly,
\begin{gather*}
\bigl( R_{\mathbf{a}} f\bigr)(x) = \frac{1}{2ix}\left( (a+ix)(b+ix)(c+ix)(d+ix) f\left(x-\frac12i\right)\right.\\
\left. \hphantom{\bigl( R_{\mathbf{a}} f\bigr)(x) =}{} - (a-ix)(b-ix)(c-ix)(d-ix) f\left(x+\frac12i\right)\right).
\end{gather*}
By \cite[equation~(1.1.10)]{KoekS} the polynomials $p_n^{(\nu)}$ of Section~\ref{sec:general-setup} correspond precisely to the Wilson polynomials.

By induction we see that the Leibniz formula for $\frac{\delta}{\delta x^2}$ is
\begin{gather}\label{eq:Wilson-Leibniz}
\left(\frac{\de}{\de x^2}\right)^n (fg)(x) = \sum_{k=0}^n \binom{n}{k} \left( \big(S^+\big)^{n-k} \left(\frac{\de}{\de x^2}\right)^{k}f\right)(x) \left( (S^-)^k \left(\frac{\de}{\de x^2}\right)^{n-k}g\right)(x),
\end{gather}
where $S^\pm f(x) = f\big(x\pm \frac12 i\big)$, assuming that the functions $f$ and $g$ are defined in a sufficiently large strip $|\Im z|< \frac12 n+\ep$. We then identify $\al^n_k$ with $\binom{n}{k}$, $\eta$ with $S^-$ and $T_{k,n} = (S^+)^{n-k} \bigl( \frac{\delta}{\delta x^2}\bigr)^k$. We check that $\frac{\eta^k (w_{\nu+k\si})(x)}{w_{\nu}(x)}$ corresponds to
\begin{gather*}
\frac{(S^-)^k \om(x;\mathbf{a}+k\si)}{\om(x;\mathbf{a})} = (-1)^k (a+ix)_k (b+ix)_k (c+ix)_k (d+ix)_k.
\end{gather*}
Then Theorem \ref{lem:gs-expansion} gives the operational expansion
\begin{gather*}
\bigl( R_{\mathbf{a}}R_{\mathbf{a}+\frac12} \cdots R_{\mathbf{a}+(n-1)\frac12}f\bigr) (x)= \sum_{k=0}^n \binom{n}{k} (-1)^k(a+ix)_k (b+ix)_k (c+ix)_k (d+ix)_k \nonumber\\
\qquad{} \times W_{n-k}\left(\left(x-\frac12 ik\right)^2; \mathbf{a}+k\frac12\right) \left(\big(S^+\big)^{n-k} \left( \frac{\delta}{\delta x^2}\right)^k f\right)(x).
\end{gather*}
Since, by \cite[equation~(1.1.8)]{KoekS}, $T_{k,n} p_{m}^{(\nu+n\si)}$ corresponds to
\begin{gather*}
(-m)_k (m+a+b+c+d+2n-1)_k W_{m-k}\left(\left(x+\frac12 i(n-k)\right)^2; \mathbf{a}+(n+k)\frac12\right),
\end{gather*}
Corollary \ref{cor:lem:gs-expansion} gives the following expansion for the Wilson polynomials:
\begin{gather}
 W_{m+n}\big(x^2;\mathbf{a}\big) = \sum_{k=0}^{n\wedge m} \frac{(-n)_k(-m)_k}{k!}(m+a+b+c+d+2n-1)_k (a+ix)_k (b+ix)_k \nonumber\\
 \hphantom{W_{m+n}\big(x^2;\mathbf{a}\big) =}{} \times (c+ix)_k (d+ix)_k W_{n-k}\left(\left(x-\frac12 ki\right)^2; \mathbf{a}+k\frac12\right)\nonumber\\
\hphantom{W_{m+n}\big(x^2;\mathbf{a}\big) =}{} \times W_{m-k}\left(\left(x+\frac12 i(n-k)\right)^2; \mathbf{a}+(n+k)\frac12\right).\label{eq:Wilson-expansion}
\end{gather}

It is also possible to consider the Leibniz formula \eqref{eq:Wilson-Leibniz} with $\eta$ identified with $S^+$, cf.\ Section~\ref{sec:backwardshift}. However, by the symmetry for the Wilson polynomials, the results are equivalent.

\section{Relation with the Toda lattice}\label{sec:Todalattice}

In Section~\ref{sec:intro} we recalled that the Toda equations~\eqref{eq:intro-Toda} are related to the modification of the orthogonality measure for orthogonal polynomials by $e^{-xt}$. Considering Corollary~\ref{cor:lem:gs-expansion2} we see that we can obtain the orthogonal polynomials for the modified measure $e^{-xt} d\mu^{(\nu)}$ assuming that (a) $\frac{\eta^k (w_{\nu+k\si})(x)}{w_{\nu}(x)}$ is a~polynomial of degree at most $k$, and (b) $f(x)=e^{-xt}$ is an eigenfunction of $T_{k,n}$ as arising in the Leibniz formula~\eqref{eq:gs-Leibniz}.

Let us first consider condition (b). The $T_{k,n}$ follow from the Leibniz formula~\eqref{eq:gs-Leibniz} for $\partial$, the operator from the raising operator for the orthogonal polynomials. Then we check the Askey-scheme and its $q$-analogue to see for which $e^{-xt}$ is an eigenfunction for these operators, and we only consider measures with infinite support. In the Askey-scheme, this only occurs if $\partial$ is the derivative $\frac{d}{dx}$, corresponding to the Hermite, Laguerre and Jacobi polynomials; the backward shift operator $\nabla$, corresponding to the Charlier and Meixner polynomials; the difference opera\-tor~$\frac{\de}{\de x}$, corresponding to the Meixner--Pollaczek and continuous Hahn polynomials.

The condition (a) is satisfied for the Hermite, Laguerre, Charlier, Meixner and Meixner--Pollaczek polynomials, as follows from the observations in Sections~\ref{sec:examples-classicalOPS},~\ref{sec:backwardshift},~\ref{sec:div-difference}. In all these cases we find an explicit expression for the polynomials with respect to $e^{-xt} d\mu^{(\nu)}(x)$ using Corollary~\ref{cor:lem:gs-expansion2}. However, in all of these cases the measure $e^{-xt} d\mu^{(\nu)}(x)$ can be identified with a related measure in the same family from the Askey scheme, i.e., one can glue $e^{-xt}$ onto the measure and obtain a~known orthogonality measure.

\begin{Proposition}\label{prop:explicit-Toda-solutions}
The following functions solve to Toda lattice equations \eqref{eq:intro-Toda}:
\begin{enumerate}\itemsep=0pt
\item[{\rm (i)}] $($Hermite$)$ $b_n(t) = -\frac12 t$, $c_n(t) = \frac12 n$ for $t\in \R$;
\item[{\rm (ii)}] $($Laguerre$)$ $b_n(t) = \frac{2n+\al+1}{1+t}$, $c_n(t) = \frac{n(n+\al)}{(1+t)^2}$ for $t>-1$;
\item[{\rm (iii)}] $($Charlier$)$ $b_n(t) = n+ae^{-t}$, $c_n(t) = nae^{-t}$ for $t\in \R$;
\item[{\rm (iv)}] $($Meixner$)$ $b_n(t) = \frac{n(ce^{-t}+1)+\be ce^{-t}}{1-ce^{-t}}$, $c_n(t) = \frac{n(n+\be-1)ce^{-t}}{ (1-ce^{-t})^2}$ for $t>\ln(c)$;
\item[{\rm (v)}] $($Meixner--Pollaczek$)$ $b_n(t) = - \frac{n+\la}{\tan(\vp-\frac12t)}$, $c_n(t) = \frac{n(n+2\la-1)}{4\sin^2(\vp-\frac12 t)}$ for $2\vp-2\pi<t<2\vp$;
\item[{\rm (vi)}] $($Krawtchouk$)$ $b_n(t) = \frac{pe^{-t}(N-n) + n(1-p)}{1+p(e^{-t}-1)}$, $c_n(t) = \frac{n(N+1-n)e^{-t}p(1-p)}{(1+p(e^{-t}-1))^2}$, $n\in \{0,1,\dots, N\}$.
\end{enumerate}
\end{Proposition}

Note that all solutions of Proposition \ref{prop:explicit-Toda-solutions} correspond to solutions of the Toda lattice equa\-tions~\eqref{eq:intro-Toda} which are separated. These solutions were obtained by Kametaka \cite{Kame,Kame-II}.

\begin{proof} The first five cases follow from the observations in Sections~\ref{sec:examples-classicalOPS}, \ref{sec:backwardshift}, \ref{sec:div-difference} and the corresponding recurrence coefficients in the relations for the monic orthogonal polynomials as in, e.g.,~\cite{KoekS}.

The final one for the Krawtchouk polynomials follows by observing that the exponential modification of Krawtchouk weight $\binom{N}{x}p^x(1-p)^{N-x}$, see, e.g., \cite[equations~(1.10.2)]{KoekS}, is of the same form with $p$ replaced by $\frac{pe^{-t}}{1+p(e^{-t}-1)}$, so that \cite[equations~(1.10.4)]{KoekS} gives the solution in the last case.
\end{proof}

In the $q$-Askey-scheme none of these operators works to give explicit orthogonal polynomials for the modification of the weight by $e^{-xt}$. In Section~\ref{sec:examples-qAskey-qderivative} we consider the case with a $q$-exponential, but this is not related to the Lax pair for the Toda lattice.

\section[Burchnall's identities for the $q$-Hahn scheme]{Burchnall's identities for the $\boldsymbol{q}$-Hahn scheme}\label{sec:examples-qAskey-qderivative}

Our next objective is to give the Burchnall type identities for orthogonal polynomials from the $q$-Askey-scheme. In this section we determine the corresponding identities in the $q$-Hahn scheme, as a subscheme of the $q$-Askey scheme, see \cite[Chapter~3]{KoornLNSFQG}. The big $q$-Jacobi polynomials are on top of the $q$-Hahn scheme, and the raising and lowering operators are given in terms of the $q$-derivative $D_q$;
\begin{gather*} 
(D_qf)(x) = \frac{f(x)-f(qx)}{x-qx},\qquad x \ne 0, \qquad (D_qf)(0) = f'(0).
\end{gather*}
The Leibniz rule for the $q$-derivative is, see \cite[Exercise~1.12(iv)]{GaspR},
\begin{gather}\label{eq:qAqd-Leibniz}
D_q^n(fg) = \sum_{k=0}^n \qbinom{n}{k}{q} (T_q^k D_q^{n-k}f) D_q^kg= \sum_{k=0}^n \qbinom{n}{k}{q} D_q^{n-k}f\big(T_q^{n-k} D_q^{k}g\big),
\end{gather}
where $(T_qf)(x)=f(xq)$ is the $q$-shift and $\qbinom{n}{k}{q}=\frac{(q;q)_n}{(q;q)_k(q;q)_{n-k}}$ is the $q$-binomial. In the context of Section~\ref{sec:general-setup} we have $\partial = D_q$, and $\cD=D_{q^{-1}}$.

\subsection[Big $q$-Jacobi polynomials]{Big $\boldsymbol{q}$-Jacobi polynomials}\label{ssec:bigqJacobipols}

The big $q$-Jacobi polynomials are defined by
\begin{gather}\label{eq:defbigqJacobipols}
P_n(x;a,b,c;q) = \rphis{3}{2}{q^{-n}, abq^{n+1}, x}{aq, cq}{q;q}
\end{gather}
and using the $q$-integral $\int_{cq}^{aq} f(x)d_qx = aq(1-q)\sum_{k=0}^\infty f\big(aq^k\big)q^k - cq(1-q)\sum_{k=0}^\infty f\big(cq^k\big)q^k$, the big $q$-Jacobi polynomials are orthogonal with respect to the weight function
\begin{gather*}
w(x; a,b,c;q) = \frac{(x/a, x/c;q)_\infty}{(x,xb/c;q)_\infty}
\end{gather*}
with respect to the discrete measure $f\mapsto \int_{cq}^{aq} f(x)w(x;a,b,cd;q)d_qx$, which corresponds to the measure $d\mu^{(\nu)}$ of Section~\ref{sec:general-setup}. In the correspondence of Section~\ref{sec:general-setup}, let $\nu$ correspond to $(a,b,c)$, and $\nu+k\si$ corresponds to $\big(q^ka,q^kb,q^kc\big)$. Let $\cV=\{(a,b,c)\,|\, 0<a<q^{-1}, 0<b<q^{-1}, c<0\}$, so that the condition on $\cV$ as in Section~\ref{sec:general-setup} is satisfied. Now $R_\nu$ corresponds to $R_{a,b,c} = M_{a,b,c}^{-1} \circ D_q \circ M_{aq,bq,cq}$:
\begin{gather*}
\bigl( R_{a,b,c} f\bigr)(x) = \frac{(1-x/aq)(1-x/cq)f(x) - (1-x)(1-xb/c)f(xq)}{(1-q)x}
\end{gather*}
and the $p_n^{(\nu)}(x)$ of Section~\ref{sec:general-setup} correspond to the big $q$-Jacobi polynomials
\begin{gather}\label{eq:bigqJacobiaspnnu}
 \frac{(aq,cq;q)_n}{(ac)^n q^{n(n+1)} (1-q)^n} P_n(x;a,b,c;q)
\end{gather}
by \cite[equation~(3.5.9--10)]{KoekS}. The adjoint $L_\nu$ of $R_\nu$ equals $D_{q^{-1}}$, since $D_{q^{-1}}$ is the formal adjoint of $\partial = D_q$ for the (unweighted) $q$-integral.

\subsubsection{First form of the Leibniz formula (\ref{eq:qAqd-Leibniz})}\label{sssec:qAqd-Leibniz1}
Put $\eta=T_q$, $\al^n_k = \qbinom{n}{k}{q}$, $T_{k,n}=D_q^k$, so that $\frac{\eta^k (w_{\nu+k\si})(x)}{w_{\nu}(x)}$ corresponds to
\begin{gather*}
\frac{w\big(xq^k; q^ka,q^kb,q^kc;q\big)}{w(x;a,b,c;q)} = (x,xb/c;q)_k.
\end{gather*}
Using \eqref{eq:bigqJacobiaspnnu}, Theorem \ref{lem:gs-expansion} gives
the operational formula
\begin{gather}
\bigl(R_{a,b,c}R_{aq,bq,cq} \cdots R_{aq^{n-1},bq^{n-1}, cq^{n-1}} f\bigr)(x) =\sum_{k=0}^n \qbinom{n}{k}{q}\frac{(x,xb/c;q)_k \big(aq^{k+1}, cq^{k+1};q\big)_{n-k}}{\big(acq^{2k}\big)^{n-k} q^{(n-k)(n-k+1)}(1-q)^{n-k}}\nonumber\\
\hphantom{\bigl(R_{a,b,c}R_{aq,bq,cq} \cdots R_{aq^{n-1},bq^{n-1}, cq^{n-1}} f\bigr)(x) =}{} \times P_{n-k}\big(xq^k;aq^k,bq^k,cq^k;q\big)\big(D_q^kf\big)(x).\label{eq:bigqJacobi-expLeibniz}
\end{gather}
Corollary \ref{cor:lem:gs-expansion} for the big $q$-Jacobi polynomials then follows using \eqref{eq:bigqJacobiaspnnu} and \cite[equation~(3.5.7)]{KoekS}. We have to take into account that \cite[equation~(3.5.7)]{KoekS} is
\begin{gather}
\bigl( D_q P_n\bigr) (\cdot;a,b,c;q) = \frac{q^{1-n}(1-q^n)\big(1-abq^{n+1}\big)} {(1-q)(1-aq)1-cq)} \bigl(T_q P_{n-1}\bigr) (\cdot;aq,bq,cq;q) \quad \Longrightarrow \nonumber \\
\bigl( D_q^k P_n\bigr) (\cdot;a,b,c;q) =(-1)^k q^{\frac12k(k+1)} \frac{\big(q^{-n}, abq^{n+1};q\big)_k}{(1-q)^n(aq,cq;q)_k} \bigl(T_q^k P_{n-k}\bigr)\big(\cdot;aq^k,bq^k,cq^k;q\big) \label{eq:bigq-Jacobi-Dqpower}
\end{gather}
taking into account the relation $D_qT_q = q T_qD_q$. After simplification, Corollary~\ref{cor:lem:gs-expansion} is
\begin{gather}
P_{n+m}(x;a,b,c;q) =\sum_{k=0}^{n\wedge m} \frac{\big(q^{-n}, q^{-m}, abq^{2n+m+1}, x,xb/c;q\big)_k }{\big(q, aq, cq, aq^{n+1},cq^{n+1};q\big)_k}(ac)^k q^{k^2 + 2k +nk}\nonumber\\
\hphantom{P_{n+m}(x;a,b,c;q) =}{} \times P_{n-k}\big(xq^k;aq^k,bq^k,cq^k;q\big) P_{m-k}\big(xq^k;aq^{n+k},bq^{n+k},cq^{n+k};q\big).\label{eq:bigqJ-expansionPm+n}
\end{gather}
Note that the left hand side is obviously symmetric in $n$ and $m$, but the right hand side is not.

Corollary \ref{cor:lem:gs-expansion2} gives
\begin{gather}
\int_{cq}^{aq} \sum_{k=0}^n \qbinom{n}{k}{q}\frac{(x,xb/c;q)_k \big(aq^{k+1}, cq^{k+1};q\big)_{n-k}}{\big(acq^{2k}\big)^{n-k} q^{(n-k)(n-k+1)}(1-q)^{n-k}} \nonumber\\
\qquad{} \times P_{n-k}\big(xq^k;aq^k,bq^k,cq^k;q\big)\big(D_q^kf\big)(x) x^p w(x;a,b,c;q) d_qx =0\label{eq:bigqJacobi-generalortho}
\end{gather}
for $p<n$. Take $f(x) = e_q(-tx) = \frac{1}{(-tx;q)_\infty}$ in \eqref{eq:bigqJacobi-generalortho} assuming $t\not\in -a^{-1}q^{-\N} \cup -c^{-1}q^{-\N}$, so that there are no poles in the support of the measures involved. Then $f$ is an eigenfunction of the $q$-derivative; $D_q^k f = \big( \frac{-t}{(1-q)}\big)^k f$, so that
\begin{gather*}
\int_{cq}^{aq} \sum_{k=0}^n \qbinom{n}{k}{q} \frac{(x,xb/c;q)_k\big(aq^{k+1}, cq^{k+1};q\big)_{n-k}}{\big(acq^{2k}\big)^{n-k} q^{(n-k)(n-k+1)}(1-q)^{n-k}}\\
\qquad{} \times P_{n-k}\big(xq^k;aq^k,bq^k,cq^k;q\big) \left( \frac{-t}{(1-q)}\right)^k \frac{1}{(-tx;q)_\infty} x^p w(x;a,b,c;q) d_qx =0
\end{gather*}
for $p<n$. Take the limit $b\downarrow 0$, i.e., specialize to the big $q$-Laguerre polynomials \cite[Section~3.11]{KoekS}. So in particular, we find an expression for the orthogonal polynomials with respect to the measure
\begin{gather*}
\int_{cq}^{aq} \frac{1}{(-tx;q)_\infty} w(x;a,0,c;q) d_qx = \int_{cq}^{aq} w(x;a,-tc,c;q) d_qx,
\end{gather*}
which is again the orthogonality measure for a big $q$-Jacobi polynomial. Hence, we find the following equality up to a constant, which is determined by considering leading coefficients or by evaluating at $x=1$, after replacing $t$ by $-b/c$:
\begin{gather}\label{eq:bigqJacobi-epansion}
\sum_{k=0}^n\frac{(q^{-n}, x;q)_k}{(q,aq,cq;q)_k} (-abq^n)^k q^{\frac12k^2 +\frac32k}P_{n-k}\big(xq^k;aq^k,0,cq^k;q\big) = P_n(x;a,b,c;q).
\end{gather}
The result \eqref{eq:bigqJacobi-epansion} seems not be easily provable using a generating function, see, e.g., \eqref{eq:Laguerre-Toda-OPs} for the classical counterpart provable as a convolution identities from generating functions.

The other natural candidate for $f$ in \eqref{eq:bigqJacobi-generalortho} is $f(x) = E_q(-xt)=(xt;q)_\infty$. Note that
\begin{gather*}
D_q^k(x\mapsto E_q(-xt)) = \left(\frac{-t}{1-q}\right)^k q^{\frac12k(k-1)} \frac{(tx;q)_\infty}{(tx;q)_k},
\end{gather*}
so that we can only get orthogonal polynomials for the measure modified by the big $q$-exponential $E_q(-xt)$ if cancellation in $\frac{(x,xb/c;q)_k}{(tx;q)_k}$ occurs. This happens precisely for $t=1$ or $t=b/c$.

In particular, taking $t=b/c$ we see that the polynomials of degree $n$ defined by
\begin{gather*}
\sum_{k=0}^n \qbinom{n}{k}{q} \frac{(x;q)_k \big(aq^{k+1}, cq^{k+1};q\big)_{n-k}}{\big(acq^{2k}\big)^{n-k} q^{(n-k)(n-k+1)}(1-q)^{n-k}} P_{n-k}\big(xq^k;aq^k,bq^k,cq^k;q\big) \left(\frac{-b/c}{1-q}\right)^k q^{\frac12k(k-1)}
\end{gather*}
are orthogonal with respect to
\begin{gather*}
f\mapsto \int_{cq}^{aq} f(x) (xb/c;q)_\infty w(x;a,b,c;q) d_qx = \int_{cq}^{aq} f(x) \frac{(x/a, x/c;q)_\infty}{(x;q)_\infty} d_qx.
\end{gather*}
Since this is the measure for the big $q$-Laguerre polynomials, i.e., big $q$-Jacobi polynomials \eqref{eq:defbigqJacobipols} with $b=0$, see, e.g., \cite[Section~4.11]{KoekS}, we also know the orthogonal polynomials as
big $q$-Laguerre polynomials. This gives the following expansion up to a constant, which is determined by evaluating at $x=1$:
\begin{gather}\label{eq:expansion-bigqLaguerre}
\sum_{k=0}^n \frac{(q^{-n}, x;q)_k}{(q, aq, cq;q)_k}(-ab)^k q^{k(k+n+1)}P_{n-k}\big(xq^k;aq^k,bq^k,cq^k;q\big) = P_n(x;a,0,c;q)
\end{gather}
and \eqref{eq:expansion-bigqLaguerre} can be considered as a kind of inverse to~\eqref{eq:bigqJacobi-epansion}.

The other case, $t=1$, can be dealt with similarly. We see that the polynomials
\begin{gather*}
\sum_{k=0}^n \qbinom{n}{k}{q} \frac{(xb/c;q)_k \big(aq^{k+1}, cq^{k+1};q\big)_{n-k}}{\big(acq^{2k}\big)^{n-k} q^{(n-k)(n-k+1)}(1-q)^{n-k}} P_{n-k}\big(xq^k;aq^k,bq^k,cq^k;q\big) \left(\frac{-1}{1-q}\right)^k q^{\frac12k(k-1)}
\end{gather*}
are orthogonal with respect to the measure
\begin{gather*}
f\mapsto \int_{cq}^{aq} f(x) (x;q)_\infty w(x;a,b,c;q) d_qx = \int_{cq}^{aq} f(x) \frac{(x/a, x/c;q)_\infty}{(xb/c;q)_\infty} d_qx \\
\hphantom{f\mapsto}{} = \left\vert \frac{c}{b}\right\vert \int_{abq/c}^{bq} f\left( \frac{xc}{b}\right) \frac{(x/b, xc/ab;q)_\infty}{(x;q)_\infty} d_qx
\end{gather*}
after a change of variable in the $q$-integral. This measure is the orthogonality measure for the big $q$-Laguerre polynomials, or, more precisely, for the big $q$-Jacobi polynomials with $(a,b,c)$ replaced by $(b,0,ab/c)$. So we find \eqref{eq:expansion-bigqLaguerre2} for the orthogonal polynomials orthogonal with respect to this measure, where the constant is determined by evaluating at $x=c/b$ and using the $q$-Saalsch\"utz formula \cite[equation~(1.7.2)]{GaspR}:
\begin{gather}
 \sum_{k=0}^n \frac{(q^{-n}, xb/c;q)_k}{q, aq, cq;q)_k} (ac)^k q^{k(k+n+1)} P_{n-k}\big(xq^k;aq^k,bq^k,cq^k;q\big) \nonumber\\
\qquad{} = \left( \frac{c}{b}\right)^n \frac{(bq, abq/c;q)_n}{(aq, cq;q)_n} P_n(xb/c;b,0,ab/c;q).\label{eq:expansion-bigqLaguerre2}
\end{gather}

\subsubsection{Second form of the Leibniz formula (\ref{eq:qAqd-Leibniz})}
As in Section~\ref{sec:backwardshift} we can consider the Leibniz expansion \eqref{eq:qAqd-Leibniz} in two different ways to fit it in the form \eqref{eq:gs-Leibniz}. Now we take $\eta=I$, $\al^n_k = \qbinom{n}{k}{q}$, $T_{k,n}=T_q^{n-k}D_q^k$, so that $\frac{\eta^k (w_{\nu+k\si})(x)}{w_{\nu}(x)}$ corresponds to
\begin{gather*}
\frac{w\big(x; q^ka,q^kb,q^kc;q\big)}{w(x;a,b,c;q)} =\big(xq^{-k}/a,xq^{-k}/c;q\big)_k
\end{gather*}
Using \eqref{eq:bigqJacobiaspnnu}, Theorem \ref{lem:gs-expansion} gives
\begin{gather*}
\bigl(R_{a,b,c}R_{aq,bq,cq} \cdots R_{aq^{n-1},bq^{n-1}, cq^{n-1}} f\bigr)(x) \nonumber\\
\qquad {} =\sum_{k=0}^n \qbinom{n}{k}{q} \frac{\big(xq^{-k}/a,xq^{-k}/c;q\big)_k \big(aq^{k+1}, cq^{k+1};q\big)_{n-k}}{\big(acq^{2k}\big)^{n-k} q^{(n-k)(n-k+1)}(1-q)^{n-k}}\nonumber\\
\qquad\quad{} \times P_{n-k}\big(x;aq^k,bq^k,cq^k;q\big)\big(T_q^{n-k}D_q^kf\big)(x)
\end{gather*}
as an alternative for \eqref{eq:bigqJacobi-expLeibniz}.

Corollary \ref{cor:lem:gs-expansion}, by \eqref{eq:bigqJacobiaspnnu} and \eqref{eq:bigq-Jacobi-Dqpower}, gives the following alternative for \eqref{eq:bigqJ-expansionPm+n};
\begin{gather*}
P_{n+m}(x;a,b,c;q) = \sum_{k=0}^{n \wedge m} \frac{\big(q^{-n}, q^{-m}, abq^{2n+m+1}, xq^{-k}/a,xq^{-k}/c;q\big)_k}{\big(q,aq,cq,aq^{n+1}, cq^{n+1};q\big)_k} (ac)^k q^{k(k+n+2)}\nonumber\\
\hphantom{P_{n+m}(x;a,b,c;q) =}{} \times P_{n-k}\big(x;aq^k,bq^k,cq^k;q\big) P_{m-k}\big(xq^n;aq^{n+k},bq^{n+k},cq^{n+k};q\big).
\end{gather*}

As in Section~\ref{sssec:qAqd-Leibniz1} we can try to relate to explicit orthogonal polynomials using Corollary~\ref{cor:lem:gs-expansion2} and a suitable choice for $f$. This orthogonality involves $\big(T_q^{n-k}D_q^kf\big)(x)$, but using $f(x)=e_q(-xt)$ or $E_q(-xt)$ does not give rise to polynomial orthogonality in
\begin{gather}
\int_{cq}^{aq} \sum_{k=0}^n \qbinom{n}{k}{q}\frac{\big(xq^{-k}/a,xq^{-k}/c;q\big)_k \big(aq^{k+1}, cq^{k+1};q\big)_{n-k}}{\big(acq^{2k}\big)^{n-k} q^{(n-k)(n-k+1)}(1-q)^{n-k}}\nonumber\\
\qquad{} \times P_{n-k}\big(x;aq^k,bq^k,cq^k;q\big)\big(T_q^{n-k}D_q^kf\big)(x) x^p w(x;a,b,c;q) d_qx =0\label{eq:bigqJacobi-generalortho2}
\end{gather}
for $p<n$, since
\begin{gather*}
\big(T_q^{n-k}D_q^k\big)e_q(-tx) = \left( \frac{-t}{1-q}\right)^k \frac{1}{(-txq^{n-k};q)_\infty} = \left( \frac{-t}{1-q}\right)^k \frac{(-tx;q)_{n-k}}{(-tx;q)_\infty}, \\
\big(T_q^{n-k}D_q^k\big)E_q(-tx) = \left( \frac{-t}{1-q}\right)^k q^{\frac12 k(k-1)} (q^ntx;q)_\infty = \left( \frac{-t}{1-q}\right)^k q^{\frac12 k(k-1)} \frac{(tx;q)_\infty}{(tx;q)_n}.
\end{gather*}
So \eqref{eq:bigqJacobi-generalortho2} cannot be rewritten as orthogonality for a family of orthogonal polynomials for these choices of $f$.

\section{Burchnall's identities for the Askey--Wilson polynomials}\label{sec:AskeyWilsoncase}

The top level in the $q$-analogue of the Askey-scheme is the family of Askey--Wilson polynomials. In order to apply the results we need the Askey--Wilson $q$-difference operator. For a func\-tion~$f$, we define $\breve{f}(z) = f\big(\frac12\big(z+z^{-1}\big)\big)$, so that $\breve{f}(z)=\breve{f}\big(z^{-1}\big)$. Then the Askey--Wilson $q$-difference operator is
\begin{gather*}
(\cD_q f)(x) = \frac{\breve{f}\big(q^{1/2}z\big)-\breve{f}\big(q^{-1/2}z\big)}{\frac12 \big(q^{1/2}-q^{-1/2}\big)\big(z-z^{-1}\big)},\qquad x= \frac12\big(z+z^{-1}\big).
\end{gather*}
Note that the denominator is the same as the numerator in case $f(x)=x$, and that the right hand side is symmetric in $z\leftrightarrow z^{-1}$. The Askey--Wilson polynomials are defined by
\begin{gather*} 
p_n(x;a,b,c,d|q) = a^{-n} (ab,ac,ad;q)_n\rphis{4}{3}{q^{-n}, abcdq^{n-1}, az, a/z}{ab, ac, ad}{q;q}, \\
 x=\frac12\big(z+z^{-1}\big).
\end{gather*}
The weight on $[-1,1]$ is given by, using $x= \cos \theta$, $z=e^{i\theta}$,
\begin{gather*}
w(x;a,b,c,d|q) = \frac{1}{\sqrt{1-x^2}} \frac{\big(z^{\pm 2};q\big)_\infty} {\big(az^{\pm 1}, bz^{\pm 1}, cz^{\pm 1}, dz^{\pm 1};q\big)_\infty} = \frac{2i}{z} \frac{\big(z^{2}, qz^{-2};q\big)_\infty} {\big(az^{\pm 1}, bz^{\pm 1}, cz^{\pm 1}, dz^{\pm 1};q\big)_\infty},
\end{gather*}
assuming that $a,b,c,d\in \R$ with $\max(|a|,|b|,|c|,|d|)<1$ and using the notation $(ez^{\pm 1};q)_\infty = (ez, e/z;q)_\infty$, see, e.g., \cite[Chapter~7]{GaspR}, \cite[Section~3.1]{KoekS}. Note that the numerator is a theta-product. So $\cV =\{ (a,b,c,d) \,|\, \max(|a|,|b|,|c|,|d|)<1 \}\subset \R^4$ in the notation of Section~\ref{sec:general-setup}. In this situation $\si$ corresponds to $\big(q^{1/2}, q^{1/2},q^{1/2},q^{1/2}\big)$ and letting $\nu$ correspond to $\mathbf{a}=(a,b,c,d)$, we identify $\nu+n\si$ with $\mathbf{a}q^{n/2}=\big(aq^{n/2},bq^{n/2},cq^{n/2},dq^{n/2}\big)$.

The raising operator $R_\nu$ corresponds to $R_{\mathbf{a}} = M^{-1}_{\mathbf{a}} \circ \cD_q \circ M_{\mathbf{a}q^{1/2}}$. Explicitly, again using $x=\frac12\big(z+z^{-1}\big)$,
\begin{gather*}
(R_{\mathbf{a}} f)(x) = A(z) \breve{f}\big(q^{1/2}z\big) - A\big(z^{-1}\big)\breve{f}\big(q^{-1/2}z\big), \\
 A(z)=\frac{-2(1-az)(1-bz)(1-cz)(1-dz)}{z(1-q)\big(1-z^2\big)}
\end{gather*}
Using \cite[equation~(3.1.12)]{KoekS} we see that $p^{(\nu)}_n$ corresponds to
\begin{gather}\label{eq:AWpols-pnnuS2}
2^n (q-1)^{-n} q^{-\frac14 n(n-1)} p_n(x; a,b,c,d\,|\, q).
\end{gather}
The corresponding Leibiz formula for $\cD_q$ is given by, see \cite[equation~(1.22)]{Isma-LAW},
\begin{gather}\label{eq:AWdiff-Leibniz}
\cD_q^n (fg) = \sum_{k=0}^n \qbinom{n}{k}{q} q^{\frac12 k(k-n)} \bigl( \eta_q^k \cD_q^{n-k} f\bigr) \bigl( \eta_q^{k-n} \cD_q^{k} g\bigr),
\end{gather}
where $(\eta_q f)(x) = \breve{f}\big(q^{1/2}z\big)$, $x=\frac12\big(z+z^{-1}\big)$. So in the context of Section~\ref{sec:general-setup} we have $\al^n_k = \qbinom{n}{k}{q} q^{\frac12k(k-n)}$, $\eta= \eta_q$ and $T_{k,n} =\eta_q^{k-n} \cD_q^k$. Then $\frac{\eta^k (w_{\nu+k\si})(x)}{w_{\nu}(x)}$ corresponds to
\begin{gather*}
(-1)^k q^{-\frac12 k^2} z^{-2k} (az, bz, cz, dz;q)_k.
\end{gather*}
Applying Theorem \ref{lem:gs-expansion} we get the following operational expansion, after simplification,
\begin{gather}
\bigl( R_{\mathbf{a}}R_{\mathbf{a}q^{1/2}} \cdots R_{\mathbf{a}q^{(n-1)/2}}f\bigr)(x)\nonumber\\
\qquad{} = \frac{2^nq^{\frac14n(n-1)}}{(q-1)^n} \sum_{k=0}^n \qbinom{n}{k}{q}\frac{q^{-nk} q^{\frac14 k(k+1)}}{2^{k}} (1-q)^{k}(az, bz, cz, dz;q)_k \label{eq:AW-operational}\\
\qquad\quad{} \times z^{-2k} p_{n-k}\left(\frac12\big(q^{k/2}z+ q^{-k/2}z^{-1}\big);aq^{k/2}, bq^{k/2}, cq^{k/2}, dq^{k/2}|q\right)\bigl( \eta_q^{k-n} \cD_q^k f\bigr)(x).\nonumber
\end{gather}

To make Corollary \ref{cor:lem:gs-expansion} explicit we iterate \cite[equation~(3.1.9)]{KoekS}
\begin{gather*}
\cD_q^k p_n (x;a,b,c,d|q) \\
\qquad{} = \frac{(-2)^k q^{\frac14 k(2n-k+3)}}{(1-q)^k} \big(q^{-m},abcdq^{n-1};q\big)_k p_{n-k}\big(x;aq^{k/2},bq^{k/2},cq^{k/2},dq^{k/2}|q\big),
\end{gather*}
so that we find the following expansion from Corollary~\ref{cor:lem:gs-expansion} for the Askey--Wilson polynomials using~\eqref{eq:AWpols-pnnuS2} and $x=\frac12\big(z+z^{-1}\big)$:
\begin{gather}
 p_{n+m}(x;a,b,c,d|q) = \sum_{k=0}^{n\wedge m}\frac{\big(q^{-n},q^{-m}, abcdq^{2n+m-1};q\big)_k}{(q;q)_k}z^{-2k}(az,bz,cz,dz;q)_k\nonumber\\
 \qquad{} \times q^{- k^2 + k +\frac12nm + \frac12 km+nk} p_{n-k}\left(\frac12\big(q^{k/2}z+q^{-k/2}z^{-1}\big);aq^{k/2},bq^{k/2},cq^{k/2},dq^{k/2}|q\right) \label{eq:AW-expansion}\\
 \qquad{} \times p_{m-k}\left(\frac12\big(q^{(k-n)/2}z+q^{(n-k)/2}z^{-1}\big);aq^{(k+n)/2},bq^{(k+n)/2},cq^{(k+n)/2},dq^{(k+n)/2}|q\right). \nonumber
\end{gather}
Note that \eqref{eq:AW-expansion} is the $q$-analogue of \eqref{eq:Wilson-expansion}. In this case we can also reinterpret the Leibniz formu\-la~\eqref{eq:AWdiff-Leibniz} as in Section~\ref{sec:backwardshift}, but the resulting identities are equivalent because of the invariance $z\leftrightarrow z^{-1}$ for the Askey--Wilson polynomials.

Since we have the Askey--Wilson polynomials on top of the $q$-analogue of the Askey-scheme, we can easily specialize \eqref{eq:AW-operational} and \eqref{eq:AW-expansion} to its subclasses and limiting cases, such as e.g. the continuous $q$-Hahn polynomials ($d=0$), the Al-Salam--Chihara polynomials ($c=d=0$), Rogers's continuous $q$-ultraspherical polynomials $(a,b,c,d)=\big(\sqrt{\be},-\sqrt{\be},q^{1/2}\sqrt{\be},-q^{1/2}\sqrt{\be}\big)$, see \cite[Section~4]{KoekS} for more limit transitions.

Note that an interesting special case is obtained by taking the case that $a=b=c=d=0$ leading to the continuous $q$-Hermite polynomials $H_n(x|q)$, which is the limit $\be\to 0$ of the continuous $q$-ultraspherical polynomials. Explicitly,
\begin{gather*}
\lim_{a\to 0} p_n\left(\frac12\big(z+z^{-1}\big);a,0,0,0|q\right) = \lim_{a\to 0} (az;q)_n z^{-n} \rphis{2}{1}{q^{-n}, 0}{q^{1-n}/az}{q, \frac{qz}{a}}\\
\hphantom{\lim_{a\to 0} p_n\left(\frac12\big(z+z^{-1}\big);a,0,0,0|q\right)}{} = z^{-n} \rphis{2}{0}{q^{-n},0}{-}{q, q^nz^2} =H_n\left(\frac12\big(z+z^{-1}\big)|q\right),
\end{gather*}
using the rewrite of the Al-Salam--Chihara polynomials as a ${}_2\vp_1$-series, cf.~\cite[Section~3.8]{KoekS}. In particular, \eqref{eq:AW-operational} then gives
\begin{gather*}
 2^{-n}(q-1)^{n}q^{-\frac14n(n-1)} \bigl( R_{\mathbf{0}}^n f\bigr)(x) \\
 = \sum_{k=0}^n \qbinom{n}{k}{q}q^{-nk} q^{\frac14 k(k+1)} 2^{-k} (1-q)^{k}z^{-2k} H_{n-k}\left(\frac12\big(q^{k/2}z+ q^{-k/2}z^{-1}\big)|q\right)\bigl( \eta_q^{k-n} \cD_q^k f\bigr)(x) \nonumber\\
=\sum_{k=0}^n \frac{(q^{-n};q)_k}{(q;q)_k}\left(-\frac12\right)^k q^{-\frac12nk} q^{\frac14 k(k+3)} (1-q)^{k}z^{-n-k} \rphis{2}{0}{q^{-n},0}{-}{q,q^nz^2}\bigl( \eta_q^{k-n} \cD_q^k f\bigr)(x).\nonumber
\end{gather*}
The special case of \eqref{eq:AW-expansion} for the continuous $q$-Hermite polynomials is
\begin{gather}
H_{n+m}(x|q) =\sum_{k=0}^{n\wedge m} \frac{(q^{-n},q^{-m};q)_k}{(q;q)_k}z^{-2k}q^{- k^2 + k +\frac12nm + \frac12 km+nk} \label{eq:contqHermite-expansion}\\
\hphantom{H_{n+m}(x|q) =}{} \times H_{n-k}\left(\frac12\big(q^{k/2}z+q^{-k/2}z^{-1}\big)|q\right) H_{m-k}\left(\frac12\big(q^{(k-n)/2}z+q^{(n-k)/2}z^{-1}\big)|q\right).\nonumber
\end{gather}
Note that \eqref{eq:contqHermite-expansion} is different from the inverse linearization formula for the continuous $q$-Hermite polynomials as derived from the $\be\to 0$ limit of the corresponding inverse linearization formula for the continuous $q$-ultraspherical polynomials in \cite[Exercise~8.12]{GaspR}.

\subsection*{Acknowledgements}
Mourad Ismail acknowledges the research support and hospitality of Radboud University during his visits which initiated this collaboration. Erik Koelink gratefully acknowledges the support of FaMAF as invited professor at Universidad Nacional de C\'ordoba for a research visit. The work of Pablo Rom\'an was supported by Radboud Excellence Fellowship, CONICET grant PIP 112-200801-01533, FONCyT grant PICT 2014-3452 and by SeCyT-UNC. We thank a~referee for useful advice.

\pdfbookmark[1]{References}{ref}
\LastPageEnding


\begin{thebibliography}{99}
\footnotesize\itemsep=0pt

\bibitem{AharAFG}
Aharmim B., Amal E.H., Fouzia E.W., Ghanmi A., Generalized {Z}ernike
 polynomials: operational formulae and generating functions, \href{https://doi.org/10.1080/10652469.2015.1012510}{\textit{Integral
 Transforms Spec. Funct.}} \textbf{26} (2015), 395--410, \href{https://arxiv.org/abs/1312.3628}{arXiv:1312.3628}.

\bibitem{AlSa}
Al-Salam W.A., Operational representations for the {L}aguerre and other
 polynomials, \href{https://doi.org/10.1215/S0012-7094-64-03113-8}{\textit{Duke Math.~J.}} \textbf{31} (1964), 127--142.

\bibitem{AndrAR}
Andrews G.E., Askey R., Roy R., Special functions, \href{https://doi.org/10.1017/CBO9781107325937}{\textit{Encyclopedia of
 Mathematics and its Applications}}, Vol.~71, Cambridge University Press,
 Cambridge, 1999.

\bibitem{Aske-SIAM}
Askey R., Orthogonal polynomials and special functions, \textit{Reg. Conf. Ser.
 Appl. Math.}, Vol.~21, Society for Industrial and Applied Mathematics,
 Philadelphia, Pa., 1975.

\bibitem{BabeBT}
Babelon O., Bernard D., Talon M., Introduction to classical integrable systems,
 \href{https://doi.org/10.1017/CBO9780511535024}{\textit{Cambridge Monographs on Mathematical Physics}}, Cambridge University Press,
 Cambridge, 2003.

\bibitem{BasoCE}
Basor E., Chen Y., Ehrhardt T., Painlev\'e~{V} and time-dependent {J}acobi
 polynomials, \href{https://doi.org/10.1088/1751-8113/43/1/015204}{\textit{J.~Phys.~A: Math. Theor.}} \textbf{43} (2010), 015204,
 25~pages, \href{https://arxiv.org/abs/0905.2620}{arXiv:0905.2620}.

\bibitem{Burc}
Burchnall J.L., A note on the polynomials of {H}ermite, \href{https://doi.org/10.1093/qmath/os-12.1.9}{\textit{Quart.~J.
 Math., Oxford Ser.}} \textbf{12} (1941), 9--11.

\bibitem{Carl}
Carlitz L., A note on the {L}aguerre polynomials, \href{https://doi.org/10.1307/mmj/1028998429}{\textit{Michigan Math.~J}}
 \textbf{7} (1960), 219--223.

\bibitem{CasaMN}
Casas F., Murua A., Nadinic M., Efficient computation of the {Z}assenhaus
 formula, \href{https://doi.org/10.1016/j.cpc.2012.06.006}{\textit{Comput. Phys. Commun.}} \textbf{183} (2012), 2386--2391,
 \href{https://arxiv.org/abs/1204.0389}{arXiv:1204.0389}.

\bibitem{GaspR}
Gasper G., Rahman M., Basic hypergeometric series, \href{https://doi.org/10.1017/CBO9780511526251}{\textit{Encyclopedia of
 Mathematics and its Applications}}, Vol.~96, 2nd~ed., Cambridge University
 Press, Cambridge, 2004.

\bibitem{Gekh}
Gekhtman M., Hamiltonian structure of non-abelian {T}oda lattice, \href{https://doi.org/10.1023/A:1007579806383}{\textit{Lett.
 Math. Phys.}} \textbf{46} (1998), 189--205.

\bibitem{GoulH}
Gould H.W., Hopper A.T., Operational formulas connected with two
 generalizations of {H}ermite polynomials, \href{https://doi.org/10.1215/S0012-7094-62-02907-1}{\textit{Duke Math.~J.}} \textbf{29}
 (1962), 51--63.

\bibitem{Isma-LAW}
Ismail M.E.H., The {A}skey--{W}ilson operator and summation theorems, in
 Mathematical Analysis, Wavelets, and Signal Processing ({C}airo, 1994),
 \href{https://doi.org/10.1090/conm/190/02300}{\textit{Contemp. Math.}}, Vol.~190, Editors M.E.H.~Ismail, M.Z.~Nashed, A.I.~Zayed, A.F.~Ghaleb, Amer. Math. Soc., Providence, RI, 1995, 171--178.

\bibitem{Isma}
Ismail M.E.H., Classical and quantum orthogonal polynomials in one variable,
 \textit{Encyclopedia of Mathematics and its Applications}, Vol.~98, Cambridge
 University Press, Cambridge, 2009.

\bibitem{IsmaKR}
Ismail M.E.H., Koelink E., Rom\'an P., {i}n preparation.

\bibitem{Kame}
Kametaka Y., On the {E}uler--{P}oisson--{D}arboux equation and the {T}oda
 equation.~{I}, \href{https://doi.org/10.3792/pjaa.60.145}{\textit{Proc. Japan Acad. Ser.~A Math. Sci.}} \textbf{60}
 (1984), 145--148.

\bibitem{Kame-II}
Kametaka Y., On the {E}uler--{P}oisson--{D}arboux equation and the {T}oda
 equation.~{II}, \href{https://doi.org/10.3792/pjaa.60.181}{\textit{Proc. Japan Acad. Ser.~A Math. Sci.}} \textbf{60}
 (1984), 181--184.

\bibitem{KoekLS}
Koekoek R., Lesky P.A., Swarttouw R.F., Hypergeometric orthogonal polynomials
 and their {$q$}-analogues, \href{https://doi.org/10.1007/978-3-642-05014-5}{\textit{Springer Monographs in Mathematics}},
 Springer-Verlag, Berlin, 2010.

\bibitem{KoekS}
Koekoek R., Swarttouw R.F., The Askey-scheme of hypergeometric orthogonal
 polynomials and its $q$-analogue, Report 98-17, Faculty of Technical
 Mathematics and Informatics, Delft University of Technology, 1998,
 \url{http://aw.twi.tudelft.nl/~koekoek/askey/}.

\bibitem{KoelVdJ}
Koelink H.T., Van Der~Jeugt J., Convolutions for orthogonal polynomials from
 {L}ie and quantum algebra representations, \href{https://doi.org/10.1137/S003614109630673X}{\textit{SIAM~J. Math. Anal.}}
 \textbf{29} (1998), 794--822, \href{https://arxiv.org/abs/q-alg/9607010}{q-alg/9607010}.

\bibitem{KoornLNSFQG}
Koornwinder T.H., Compact quantum groups and {$q$}-special functions, in
 Representations of {L}ie Groups and Quantum Groups ({T}rento, 1993),
 \textit{Pitman Res. Notes Math. Ser.}, Vol.~311, Longman Sci. Tech., Harlow,
 1994, 46--128.

\bibitem{Magn}
Magnus W., On the exponential solution of differential equations for a linear
 operator, \href{https://doi.org/10.1002/cpa.3160070404}{\textit{Comm. Pure Appl. Math.}} \textbf{7} (1954), 649--673.

\bibitem{Niel}
Nielsen N., Recherches sur les polynomes d'{H}ermite, \textit{Det Kgl. Danske
 Videnskabernes Selskab. Math.-Fys. Meddelelser~I} \textbf{6} (1918), 1--78.

\bibitem{Rain}
Rainville E.D., Special functions, The Macmillan Co., New York, 1960.

\bibitem{Sing}
Singh R.P., Operational formulae for {J}acobi and other polynomials,
 \textit{Rend. Sem. Mat. Univ. Padova} \textbf{35} (1965), 237--244.

\bibitem{Zhed}
Zhedanov A., The Toda chain: solutions with dynamical symmetry and classical
 orthogonal polynomials, \href{https://doi.org/10.1007/BF01028245}{\textit{Theoret. and Math. Phys.}} \textbf{82} (1990),
 6--11.

\bibitem{Zhed-R}
Zhedanov A., Elliptic solutions of the {T}oda chain and a generalization of the
 {S}tieltjes--{C}arlitz polynomials, \href{https://doi.org/10.1007/s11139-013-9515-x}{\textit{Ramanujan~J.}} \textbf{33} (2014),
 157--195, \href{https://arxiv.org/abs/0712.0058}{arXiv:0712.0058}.

\end{thebibliography}
\end{document}